\documentclass[12pt,reqno,a4paper]{amsart}
\usepackage[a4paper,margin=1in]{geometry}

\title[Debarre--Voisin varieties and Mukai homomorphisms]{K3 Categories and Mukai Homomorphisms\\for Debarre--Voisin Varieties}
\author{Ryo Akazawa}

\address{
Department of Mathematics, Faculty of Science and Engineering,
Waseda University,
3-4-1 Ohkubo, Shinjuku, Tokyo 169-8555, Japan}
\email{rakazawa11235@ruri.waseda.jp}

\usepackage[T1]{fontenc}
\usepackage{amsmath,amssymb}
\usepackage{mathtools}
\usepackage{newtxtext}
\usepackage[nosymbolsc]{newtxmath}
\usepackage[mathscr]{eucal}

\usepackage{graphicx}
\usepackage{array,booktabs}
\usepackage{listings}
\usepackage{tikz-cd}

\usepackage{xcolor}
\usepackage{hyperref}
\hypersetup{%
 colorlinks=true,%
 linkcolor=blue!75!black,%
 citecolor=red!70!black,%
 urlcolor=teal!85!black%
}
\usepackage[initials,alphabetic,lite]{amsrefs}
\makeatletter
\def\bibtex@style{amsra-shorthand}
\makeatother
\theoremstyle{plain}
\newtheorem{thm}{Theorem}[section]
\newtheorem{lem}[thm]{Lemma}
\newtheorem{prop}[thm]{Proposition}
\newtheorem{cor}[thm]{Corollary}
\newtheorem{conj}[thm]{Conjecture}

\newtheorem*{thm*}{Theorem}
\newtheorem*{cor*}{Corollary}

\theoremstyle{definition}
\newtheorem{dfn}[thm]{Definition}

\theoremstyle{remark}

\newtheorem{rem}[thm]{Remark}

\DeclareMathOperator{\Ku}{Ku}

\DeclareMathOperator{\Hom}{Hom}

\DeclareMathOperator{\Ext}{Ext}
\DeclareMathOperator{\RHom}{RHom}
\DeclareMathOperator{\ch}{ch}
\DeclareMathOperator{\td}{td}

\DeclareMathOperator{\Gr}{Gr}

\newcommand{\Db}{\mathrm{D}^{\mathrm{b}}}

\newcommand{\C}{\mathbb{C}}

\newcommand{\Q}{\mathbb{Q}}
\newcommand{\Z}{\mathbb{Z}}

\newcommand{\Ktop}{K_{\mathrm{top}}}
\newcommand{\HH}{\mathrm{HH}}
\newcommand{\cE}{\mathcal E}
\newcommand{\cF}{\mathcal F}
\newcommand{\cP}{\mathcal P}
\newcommand{\cZ}{\mathcal Z}
\newcommand{\uDV}{u_{\mathrm{DV}}}
\newcommand{\Hyp}{\mathbf U}

\newcommand{\NHdg}{N_{\mathrm{Hdg}}}

\makeatletter
\def\@secnumfont@section{}
\def\@secnumfont@subsection{\protect\bfseries}
\def\@seccntformat#1{%
  \protect\textup{\protect\@secnumfont
    \csname @secnumfont@#1\endcsname
    \csname the#1\endcsname
    \protect\@secnumpunct
  }%
}
\makeatother

\numberwithin{equation}{section}

\providecommand{\Q}{\mathbb Q}
\providecommand{\C}{\mathbb C}
\providecommand{\Z}{\mathbb Z}

\providecommand{\Gr}{\operatorname{Gr}}

\providecommand{\cE}{\mathcal E}

\providecommand{\ch}{\operatorname{ch}}
\providecommand{\Hom}{\operatorname{Hom}}

\BibSpec{collection.article}{%
    +{}  {\PrintAuthors}                {author}
    +{,} { \textit}                     {title}
    +{.} { }                            {part}
    +{:} { \textit}                     {subtitle}
    +{,} { \PrintContributions}         {contribution}
    +{,} { \PrintConference}            {conference}
    +{}  {\PrintBook}                   {book}
    +{,} { }                            {booktitle}
    +{,} { }                            {series}
    +{,} { \voltext}                    {volume}
    +{,} { }                            {publisher}
    +{,} { \PrintDateB}                 {date}
    +{,} { pp.~}                        {pages}
    +{,} { }                            {status}
    +{,} { \PrintDOI}                   {doi}
    +{,} { available at \eprint}        {eprint}
    +{}  { \parenthesize}               {language}
    +{}  { \PrintTranslation}           {translation}
    +{;} { \PrintReprint}               {reprint}
    +{.} { }                            {note}
    +{.} {}                             {transition}
    +{}  {\SentenceSpace \PrintReviews} {review}
}

\providecommand{\Db}{\mathrm{D}^{\mathrm b}}
\providecommand{\Ku}{\operatorname{Ku}}
\providecommand{\td}{\operatorname{td}}
\providecommand{\Ext}{\operatorname{Ext}}
\providecommand{\RHom}{\operatorname{RHom}}
\providecommand{\Ktop}{K_{\mathrm{top}}}
\providecommand{\HH}{\mathrm{HH}}
\providecommand{\cF}{\mathcal F}
\providecommand{\cP}{\mathcal P}
\providecommand{\cZ}{\mathcal Z}
\providecommand{\uDV}{u_{\mathrm{DV}}}
\providecommand{\Hyp}{\mathbf U}

\providecommand{\NHdg}{N_{\mathrm{Hdg}}}

\usepackage{upquote}
\renewcommand{\lstlistingname}{Code}
\providecommand{\colonequals}{\coloneqq}

\begin{document}

\begin{abstract}

Let \(M=X_6^\sigma\) be a smooth Debarre--Voisin hyperk\"ahler fourfold and
\(X=X_3^\sigma\) the associated smooth hyperplane section of \(\operatorname{Gr}(3,10)\).
Using the incidence correspondence between \(X\) and \(M\), we construct
a family \(\mathcal E\) of objects in the Kuznetsov component
\(\operatorname{Ku}(X)\). We write \(v\) for the common Mukai vector
of its fibers. We prove that the associated Mukai homomorphism
\[
 \Theta_{\mathcal E}|_{v^\perp}\colon v^\perp\xrightarrow{\sim}H^2(M,\mathbb Z)
\]
is an integral Hodge isometry and explicitly determine the inverse image
of the Pl\"ucker polarization.
For a very general point of the trivector divisor corresponding to
discriminant \(24\), we also construct an integral Hodge isometry between
the full Mukai lattices of \(\operatorname{Ku}(X)\) and the associated twisted K3 surface.

\end{abstract}

\maketitle
\setcounter{tocdepth}{1}
\tableofcontents

\section{Introduction}

\subsection{Geometric setting and main results}\label{subsec:geometric-setting}

Let \(V\) be a ten-dimensional complex vector space. For
\(0\ne\sigma\in\bigwedge^3V^\vee\), set
\[
 X=X_3^\sigma\colonequals \{W\in\Gr(3,V)\mid\sigma|_W=0\},\qquad
 M=X_6^\sigma\colonequals \{U\in\Gr(6,V)\mid\sigma|_U=0\}.
\]
Throughout, we assume that both \(X\) and \(M\) are smooth. Then \(M\) is a
four-dimensional \emph{irreducible holomorphic symplectic variety} of
\(\mathrm{K3}^{[2]}\)-type, and its Pl\"ucker class
\(h=c_1(\mathcal O_M(1))\) satisfies \(q_M(h)=22\), where \(q_M\) is the
Beauville--Bogomolov form~\cite{DebarreVoisin2010Grassmann}.

The derived category of \(X\) admits a semiorthogonal decomposition
\[
 \Db(X)=\langle\Ku(X),A_1,\ldots,A_{108}\rangle.
\]
The \emph{Kuznetsov component} \(\Ku(X)\) is an example of a
\emph{K3 category}~\cites{Fonarev2013Lefschetz,Kuznetsov2019FractionalCY}.
We denote the projection functor by \(\delta\colon\Db(X)\to\Ku(X)\).

For each \(U\in M\), let \(Z_U\colonequals\Gr(3,U)\subset X\), and let
\(i_U\colon Z_U\hookrightarrow X\) be the inclusion.
Adjunction gives \(\det N_{Z_U/X}\simeq\mathcal O_{Z_U}(3)\).
Using this line bundle, we define
\[
 F_U\colonequals i_{U*}\mathcal O_{Z_U}(3),\qquad
 E_U\colonequals\delta(F_U)\in\Ku(X).
\]
This construction is analogous to the construction of
\(\delta_Y(\mathcal O_\ell(1))\) from a line \(\ell\) on a cubic fourfold
\(Y\), where \(\delta_Y\colon\Db(Y)\to\Ku(Y)\) is the standard projection.
The class of this object is \(\lambda_1\), introduced by Addington--Thomas
~\cite{AddingtonThomas2014Hodge}*{Section~2.4}.

We define the \emph{Mukai lattice} of \(\Ku(X)\) and its pairing by
\[
 L\colonequals\widetilde H(\Ku(X),\Z)
 =\{x\in K^0_{\mathrm{top}}(X)\mid\chi_X(A_i,x)=0\text{ for all }i\},
 \qquad (x,y)\colonequals-\chi_X(x,y).
\]
The topological Chern character equips \(L\) with a weight-two Hodge
structure; see Section~\ref{sec:preliminaries}.
The class of \(E_U\) is independent of \(U\), and we write
\(v\colonequals[E_U]\in L\).

Let
\[
 \cZ\colonequals\{(W,U)\in X\times M\mid W\subset U\}
\]
be the incidence variety. Denote the projections from \(X\times M\) to
\(X\) and \(M\) by \(p,q\), respectively, and the closed immersion by
\(i\colon\cZ\hookrightarrow X\times M\).
Set \(\mathcal O_{\cZ}(3)\colonequals i^*p^*\mathcal O_X(3)\), and define
\[
 \cF\colonequals i_*\mathcal O_{\cZ}(3),\qquad
 \cE\colonequals\delta_M(\cF).
\]
Here \(\delta_M\) is the relative projection functor over \(M\), and the
derived fiber of \(\cE\) at \(U\) is \(E_U\).
Using the global \(K\)-class of this family, we define, for \(x\in L\),
\[
 \Theta_{\cE}(x)\colonequals
 c_1\!\left(q_!([\cE]^\vee\otimes p^*x)\right)\in H^2(M,\Z).
\]
The symbol \(\vee\) denotes the derived dual, and \(q_!\) is the pushforward
in topological \(K\)-theory. The family and this map are constructed in
Section~\ref{sec:incidence-kernel}.

Our main result is the following.
\begin{theorem}\label{thm:main}
Let \(0\ne\sigma\in\bigwedge^3V^\vee\), and assume that
\(X=X_3^\sigma\) and \(M=X_6^\sigma\) are smooth.
Then \(v^2=2\), and the restriction of \(\Theta_{\cE}\) to \(v^\perp\),
\[
 \Theta_{\cE}|_{v^\perp}\colon v^\perp\longrightarrow H^2(M,\Z),
\]
is an integral Hodge isometry for the Mukai pairing and the
Beauville--Bogomolov form. Moreover, one can explicitly construct a
primitive class \(\uDV\in v^\perp\) such that
\[
 \uDV^2=22,\qquad \operatorname{div}_{v^\perp}(\uDV)=2,\qquad
 \Theta_{\cE}(\uDV)=h.
\]
\end{theorem}

Our second main result compares the Mukai lattice with that of a twisted
K3 surface in the discriminant-\(24\) case. In the moduli space of trivectors,
let \(\mathcal D^{1,6,10}\) denote the divisor of points admitting a flag
\(V_1\subset V_6\subset V\), with \(\dim V_i=i\), such that
\[
 \sigma(V_1,V_6,V)=0.
\]
Benedetti--Song show that this divisor maps birationally to the Heegner
divisor \(\mathcal D_{24}\) in the period space of polarized hyperk\"ahler
fourfolds~\cite{BenedettiSong2021Divisors}*{Proposition~5.2}.

For a very general point \([\sigma]\in\mathcal D^{1,6,10}\), there is a
twisted K3 surface \((S,\beta)\) associated with \(M\).
Let \(\ell\) be its degree-six polarization, and choose a \(B\)-field lift
\(B_S\) as in Section~\ref{sec:D24}.
Fix the isomorphism of Benedetti--Song
\[
 \iota\colon M\xrightarrow{\sim}M_{S,B_S,\ell}(w)
\]
~\cite{BenedettiSong2021Divisors}*{Theorem~5.19}.
The right-hand side is the moduli space of \(\beta\)-twisted
\(\ell\)-stable sheaves on \(S\) with common Mukai vector
\[
 w\colonequals(2,2B_S,0)\in\widetilde H(S,B_S,\Z).
\]
We write \(\theta_S\colon w^\perp\to H^2(M,\Z)\) for the Mukai
homomorphism defined using this isomorphism.

\begin{theorem}
\label{thm:D24-main}
Let \([\sigma]\) be a very general point of \(\mathcal D^{1,6,10}\).
Then there exists an integral Hodge isometry
\[
 \varphi\colon\widetilde H(\Ku(X),\Z)
 \xrightarrow{\sim}\widetilde H(S,B_S,\Z)
\]
such that \(\varphi(v)=\pm w\). Moreover, \(\varphi\) can be chosen so that
\[
 \varphi|_{v^\perp}=\theta_S^{-1}\circ\Theta_{\cE}|_{v^\perp}.
\]
\end{theorem}

The first theorem shows that the map associated with the relative projection
of the incidence family is an integral Hodge isometry between \(v^\perp\)
and the whole of \(H^2(M,\Z)\). In the second theorem, at a very general
point of \(\mathcal D^{1,6,10}\), we compare this map with the Mukai
homomorphism of a twisted K3 surface and extend the resulting isomorphism
\(v^\perp\simeq w^\perp\) to the full Mukai lattices on both sides.
A Fourier--Mukai equivalence induces an integral Hodge isometry
~\cite{Perry2022IntegralHodge}*{Section~5.1, Definition~6.4}.
Our result therefore provides Hodge-theoretic evidence for an equivalence
\(\Ku(X)\simeq\Db(S,\beta)\).

\subsection{Mukai theory and related results}

The Mukai lattice of a K3 surface \(S\) is
\[
 \widetilde H(S,\Z)=H^0(S,\Z)\oplus H^2(S,\Z)\oplus H^4(S,\Z),\qquad
 ((r,c,a),(r',c',a'))=c\cdot c'-ra'-r'a.
\]
The Mukai vector of an object \(E\in\Db(S)\) is
\[
 v_S(E)=\ch(E)\sqrt{\td(S)}
 =\left(r,c_1(E),r+\frac{c_1(E)^2}{2}-c_2(E)\right),
\]
and Hirzebruch--Riemann--Roch gives
\(\chi(E,F)=-(v_S(E),v_S(F))\).
For example, the ideal sheaf of a zero-dimensional subscheme
\(\xi\subset S\) of length two satisfies
\(v_S(I_\xi)=(1,0,-1)\) and
\(v_S(I_\xi)^2=2\).

Mukai constructed symplectic forms on moduli spaces of stable sheaves on K3 surfaces
~\cite{Mukai1984Symplectic}.
The relation between their second cohomology and the Mukai lattice was studied
by O'Grady and Yoshioka
~\cites{OGrady1997WeightTwo,Yoshioka2006TwistedSheaves}.

Let \(S\) be a complex projective K3 surface, let \(a\) be a primitive
algebraic Mukai vector with \(a^2>0\), and let
\(\tau\) be an \(a\)-generic Bridgeland stability condition in the distinguished
component \(\operatorname{Stab}^{\dagger}(S)\). Here genericity means that there are no strictly
semistable objects of class \(a\). If the moduli space
\(M_a=M_\tau(a)\) of stable objects is nonempty, it is a hyperk\"ahler
variety of dimension \(a^2+2\)
~\cite{BayerMacri2014MMP}*{Theorem~3.6(a)}.
The \emph{Mukai homomorphism} associated with a universal or
quasi-universal family gives an integral Hodge isometry
\[
 \theta_a\colon a^\perp\xrightarrow{\sim}H^2(M_a,\Z)
\]
~\cite{BayerMacri2014MMP}*{Theorem~3.6(b), p.~525}.
When a universal family exists, an algebraic class
\(x\in a^\perp\cap\widetilde H_{\mathrm{alg}}(S,\Z)\) and this family
define a determinant line bundle \(\lambda_x\) satisfying
\(\theta_a(x)=c_1(\lambda_x)\).

For an ample line bundle \(H\) on \(M_a\), the class
\[
 u_H\colonequals\theta_a^{-1}(c_1(H))\in a^\perp
\]
is uniquely determined. It is the inverse image of the polarization
\(c_1(H)\) in the Mukai lattice and satisfies
\(u_H^2=q_{M_a}(c_1(H))\), since \(\theta_a\) is an isometry.

In the Kuznetsov component of a smooth cubic fourfold \(Y\), a line
\(\ell\subset Y\) defines two classes
\[
 \lambda_1=[\delta_Y(\mathcal O_\ell(1))],\qquad
 \lambda_2=[\delta_Y(\mathcal O_\ell(2))].
\]
For our pairing \((\, ,\,)=-\chi\), they satisfy
\(\lambda_1^2=\lambda_2^2=2\) and \((\lambda_1,\lambda_2)=-1\), and generate
the standard \(A_2\) sublattice of the Mukai lattice
~\cite{AddingtonThomas2014Hodge}*{Section~2.4}.
If a primitive algebraic Mukai vector
\(b\in\widetilde H_{\mathrm{alg}}(\Ku(Y),\Z)\) satisfies \(b^2>0\) and
\(\tau\) is a \(b\)-generic stability condition in the distinguished component
\(\operatorname{Stab}^{\dagger}(\Ku(Y))\),
the Mukai homomorphism associated with the moduli space of stable objects,
\[
 \theta_b\colon b^\perp\xrightarrow{\sim}
 H^2\!\left(M_\tau(\Ku(Y),b),\Z\right),
\]
is an integral Hodge isometry
~\cite{PMIHES_2021__133__157_0}*{Theorem~29.2(2)}.
Here the orthogonal complement is taken in
\(\widetilde H(\Ku(Y),\Z)\).
The map \(\Theta_{\cE}|_{v^\perp}\) in Theorem~\ref{thm:main} is an analogue
of this Mukai homomorphism in the Debarre--Voisin setting.

In the Debarre--Voisin case, let \(p_{\cZ},q_{\cZ}\) be the projections
from \(\cZ\) to \(X,M\), respectively.
The incidence isometry of Benedetti--Song gives
\[
 H^{20}(X,\Z)_{\mathrm{van}}(9)
 \xrightarrow{\ q_{\cZ*}p_{\cZ}^*\ } H^2(M,\Z)_{\mathrm{prim}}
\]
~\cites{DebarreVoisin2010Grassmann,BenedettiSong2021Divisors}.
Here \(H^2(M,\Z)_{\mathrm{prim}}\colonequals h^\perp\), and \((9)\) is the
Tate twist taking weight \(20\) to weight \(2\).
The pairings are the intersection form on the source and \(-q_M\) on the
target. This correspondence identifies the rank-\(22\) vanishing and
primitive parts. In the proof of Theorem~\ref{thm:main}, we show that the
restriction of \(\Theta_{\cE}\) to the vanishing part equals the negative
of this correspondence. We then prove \(\Theta_{\cE}(\uDV)=h\), obtaining
a Hodge isometry on the full rank-\(23\) lattices after tensoring with
\(\Q\). Surjectivity over \(\Z\) follows from integrality of the map and
the discriminants of the two lattices.
To the author's knowledge, realizing the Hodge isometry by applying the relative
projection functor to the incidence family and explicitly determining the inverse
image of the Pl\"ucker polarization are new. In the discriminant-\(24\) case, we compare this map
with the Mukai homomorphism associated with the twisted K3 moduli
description of Benedetti--Song, through their common target \(H^2(M,\Z)\).

For hypersurfaces, Kuznetsov--Manivel--Markushevich studied the relation between
Fourier--Mukai transforms and Abel--Jacobi maps
~\cite{KuznetsovManivelMarkushevich2010AbelJacobi}.
Perry constructed Mukai Hodge structures for CY2 categories
~\cite{Perry2022IntegralHodge}*{Example~6.9}.
The present work continues these approaches.

\subsection{Outline of the paper}

Section~\ref{sec:preliminaries} recalls the Mukai lattice and projection
functors. In Section~\ref{sec:euler-square}, we construct \(E_U\) and
compute \(v_U^2=2\) for its class \(v_U=[E_U]\).
Section~\ref{sec:incidence-kernel} constructs the incidence family and its
relative projection \(\cE\), proves that \(v_U\) is constant, denotes its
common value by \(v\), and defines the Mukai homomorphism \(\Theta_{\cE}\).
In Section~\ref{sec:vanishing}, we show that its restriction to the
vanishing part is the negative of the incidence map.
Section~\ref{sec:algebraic-lattice} determines the inverse image \(\uDV\)
of the polarization. These two computations and integrality of the map
give the integral isometry of the main theorem in
Section~\ref{sec:main-proof}. Section~\ref{sec:D24} compares the Mukai lattices
of \(\Ku(X)\) and the associated twisted K3 surface.
Appendix~\ref{sec:computation} contains the formulas,
code, and output used in the computations.

\subsection*{Acknowledgements}
The author thanks his advisor, Professor Yasunari Nagai, for helpful
discussions, advice, and encouragement.

\section{Preliminaries}\label{sec:preliminaries}

\subsection{Notation and geometric setting}\label{subsec:notation-geometry}

All varieties and categories are defined over \(\C\).
For a smooth projective variety \(Y\), we write
\(\Db(Y)=\operatorname{Perf}(Y)\) for the bounded derived category of coherent sheaves on \(Y\).
Pullbacks, tensor products, duals, and pushforwards are derived unless otherwise stated.

Let \(K^0_{\mathrm{top}}(Y)\) denote the Grothendieck group of topological complex vector bundles.
The class of a perfect complex is defined by the alternating sum of a locally free resolution
and is additive in triangles. For \(\pi_Y\colon Y\to\operatorname{pt}\), we define the topological Euler pairing by
\[
 \chi_Y(x,y)\colonequals (\pi_Y)_!(x^\vee\otimes y).
\]
For algebraic objects, we have \(\chi_Y([E],[F])=\sum_k(-1)^k\dim\Ext_Y^k(E,F)\)
~\cite{Perry2022IntegralHodge}*{Lemma~5.2}.

We denote the tautological subbundle on \(\Gr(k,V)\) by \(\mathcal U\),
the Pl\"ucker line bundle by \(\mathcal O(1)=\det\mathcal U^\vee\),
and write \(G(t)=G\otimes\mathcal O(t)\).
The notation \(\Sigma^\lambda\) denotes the Schur functor associated with a partition \(\lambda\).
Let \(X=X_3^\sigma\) and \(M=X_6^\sigma\) be as in Section~\ref{subsec:geometric-setting},
and assume that both are smooth. We write \(\mathcal U_X\) and \(\mathcal U_M\)
for the restrictions of the tautological subbundles to \(X\) and \(M\), respectively.
Their ranks are \(3,6\), respectively.

A lattice \(J\) is a finitely generated free \(\Z\)-module equipped with a
nondegenerate symmetric bilinear form with values in \(\Z\).
We set \(J^\vee=\{x\in J_\Q\mid(x,J)\subset\Z\}\),
\(A_J=J^\vee/J\), \(|\operatorname{disc}J|=|A_J|\),
and \(\operatorname{div}_J(a)=\gcd\{(a,x)\mid x\in J\}\).
A vector \(a\) is primitive if it cannot be written as \(rb\), \(b\in J,\ r\in\Z,\ r>1\),
and \(J\) is unimodular if \(J=J^\vee\).
A ``very general point'' is a point outside a countable union of proper Zariski closed subsets.

By adjunction, \(\omega_X\simeq\mathcal O_X(-9)\).
The natural polarization \(h\) on \(M\) satisfies
\[
 q_M(h)=22,\qquad\operatorname{div}_{H^2(M,\Z)}(h)=2.
\]
The second cohomology lattice is
\[
 (H^2(M,\Z),q_M)
 \simeq\Hyp^{\oplus3}\oplus E_8(-1)^{\oplus2}\oplus\langle-2\rangle
\]
~\cites{DebarreVoisin2010Grassmann,BenedettiSong2021Divisors}.
Here \(\Hyp\) is the hyperbolic plane with Gram matrix
\(\left(\begin{smallmatrix}0&1\\1&0\end{smallmatrix}\right)\).
In particular, \(|\operatorname{disc}H^2(M,\Z)|=2\).

\subsection{Semiorthogonal decompositions and mutations}
For details on the material in this subsection, see, for example, Macr\`i--Stellari
~\cite{MacriStellari2018NoncommutativeK3}*{Section~2.1}.
In the notation \(\mathcal D=\langle\mathcal C_1,\ldots,\mathcal C_m\rangle\),
our convention is that \(\RHom(\mathcal C_j,\mathcal C_i)=0\) for \(j>i\).
A subcategory is called \emph{admissible} if its inclusion functor has both left and right adjoints,
and an object is called \emph{exceptional} if \(\RHom(A,A)\simeq\C[0]\).
A semiorthogonal decomposition
\(\Db(Y)=\langle\mathcal B,\mathcal B(1),\ldots,\mathcal B(m-1)\rangle\)
consisting of twists of the same admissible subcategory \(\mathcal B\) is called a
\emph{rectangular Lefschetz decomposition}
~\cite{MacriStellari2018NoncommutativeK3}*{Definition~2.34}.

Applying Kuznetsov's hyperplane section theorem to Fonarev's rectangular Lefschetz decomposition gives
\begin{equation}\label{eq:sod}
 \Db(X)=\langle\Ku(X),\mathcal B_X,\mathcal B_X(1),\ldots,\mathcal B_X(8)\rangle
 =\langle\Ku(X),A_1,\ldots,A_{108}\rangle
\end{equation}
~\cite{Fonarev2013Lefschetz}*{Theorem~4.1, Remark~4.7},
~\cite{Kuznetsov2019FractionalCY}*{Corollary~4.4}.
Here \(A_{\lambda,t}=\Sigma^\lambda\mathcal U_X^\vee(t)\), \(0\le t\le8\),
\[
 \Lambda=\{(0),(1),(1,1),(2),(2,1),(2,2),
              (3),(3,1),(3,2),(4),(4,1),(4,2)\},
\]
and \(\mathcal B_X=\langle A_{\lambda,0}\mid\lambda\in\Lambda\rangle\).
We enumerate these 108 exceptional objects as \(A_i\) in increasing order of twist,
using the displayed order of \(\Lambda\) within each block.
Setting \(\mathcal A_X\colonequals \langle A_1,\ldots,A_{108}\rangle\), we have
\(\Ku(X)=\mathcal A_X^\perp=\{E\mid\RHom(A_i,E)=0\text{ for all }i\}\).

For an exceptional object \(A_j\), we define the \emph{left mutation} through \(A_j\) by the evaluation triangle
\begin{equation}\label{eq:absolute-mutation}
 \RHom_X(A_j,G)\otimes A_j\longrightarrow G\longrightarrow\mathbb L_{A_j}(G)
\end{equation}
and set \(\delta\colonequals \mathbb L_{A_1}\cdots\mathbb L_{A_{108}}\).
By semiorthogonality, orthogonality to the objects through which we have already mutated is preserved
at subsequent stages, and \(\delta\) is left adjoint to the inclusion \(\Ku(X)\hookrightarrow\Db(X)\).
The projection triangle
\begin{equation}\label{eq:absolute-projection-triangle}
 P_G\longrightarrow G\longrightarrow\delta(G),\qquad P_G\in\mathcal A_X
\end{equation}
is functorially defined
~\cite{MacriStellari2018NoncommutativeK3}*{Remark~2.2, Proposition~2.6}.

A smooth proper category whose Serre functor is \([2]\) is called a \emph{Calabi--Yau \(2\) category}.
Such a category is called a \emph{K3 category} if its Hochschild homology satisfies
\(\HH_{\pm2}=\C,\ \HH_0=\C^{22}\) and vanishes in all other degrees.
It is known that \(\Ku(X)\) satisfies these conditions
~\cite{MacriStellari2018NoncommutativeK3}*{Section~3.3}.

\subsection{The Mukai lattice and its Hodge structure}

\begin{defn}\label{def:ktop-ku}
The \emph{topological K-theory} of \( \Ku(X) \) and its pairing are defined by
\[
 \Ktop(\Ku(X))\colonequals \{x\in K^0_{\mathrm{top}}(X)\mid
           \chi_X(A_i,x)=0\text{ for all }i\},\qquad
 (x,y)\colonequals -\chi_X(x,y).
\]
\end{defn}
The inclusion functor \(\Ku(X)\hookrightarrow\Db(X)\) induces an injective map
\[
 K^{\mathrm{top}}_0(\Ku(X))\longrightarrow K^0_{\mathrm{top}}(X)
\]
whose image is the group in Definition~\ref{def:ktop-ku}.
Indeed, additivity with respect to semiorthogonal decompositions gives
\(K^0_{\mathrm{top}}(X)=K^{\mathrm{top}}_0(\Ku(X))\oplus\bigoplus_i\Z[A_i]\)
~\cite{Perry2022IntegralHodge}*{Theorem~5.1, Lemma~5.2}.
Since the exceptional Euler matrix is triangular with all diagonal entries equal to one,
the orthogonality conditions in the definition are equivalent to the vanishing of all coefficients of the exceptional components.

The Mukai vector \(v_X(x)=\ch(x)\sqrt{\td(X)}\) induces an isomorphism
\(K^0_{\mathrm{top}}(X)_\C\simeq H^{2\ast}(X,\C)\).
The only off-diagonal Hodge components of \(X\) are \(H^{11,9},H^{9,11}\).
\begin{defn}\label{def:mukai-lattice}
We equip \(\Ktop(\Ku(X))\) with the following weight-two Hodge decomposition
and write \(L=\widetilde H(\Ku(X),\Z)\):
\[
 \begin{aligned}
 L^{2,0}&=v_X^{-1}(H^{11,9}(X)),\\
 L^{1,1}&=v_X^{-1}\left(\bigoplus_kH^{k,k}(X)\right)\cap L_\C,\\
 L^{0,2}&=v_X^{-1}(H^{9,11}(X)).
 \end{aligned}
\]
The pairing is that of Definition~\ref{def:ktop-ku}, and the algebraic Mukai lattice is defined by
\[
 \widetilde H_{\mathrm{alg}}(\Ku(X),\Z)
 \colonequals \operatorname{im}\left(K_0(\Ku(X))\to\Ktop(\Ku(X))\right).
\]
\end{defn}

The description in terms of topological \(K\)-theory follows Addington--Thomas,
while the sign of the Mukai pairing \((\, ,\,)=-\chi\) follows the convention of Macr\`i--Stellari
~\cites{AddingtonThomas2014Hodge,MacriStellari2018NoncommutativeK3}.
This Hodge structure is obtained by applying Perry's construction
~\cite{Perry2022IntegralHodge}*{Proposition~5.4, Definition~6.4}
to our variety \(X\).
Indeed, we restrict the isomorphism of weight-zero Hodge structures
\(K^{\mathrm{top}}_0(\Db(X))_\Q\simeq\bigoplus_kH^{2k}(X,\Q)(k)\)
to the semiorthogonal component and take the Tate twist by \((-1)\).
Since each component of \(\sqrt{\td(X)}\) is of diagonal type, the Mukai vector gives the same Hodge types.
Moreover, Serre duality and adjunction on topological \(K\)-theory
~\cite{AddingtonThomas2014Hodge}*{Section~2.1, equation~(2.4) and the following paragraph}
give
\(\chi(x,y)=\chi(y,S_{\Ku(X)}x)=\chi(y,x)\), so the pairing is symmetric.
By HRR, the pairing of classes of types \((a,b)\) and \((c,d)\) can be nonzero only if
\(a+c=b+d=2\). Thus the pairing is compatible with this Hodge structure.
The class of an object has type \((1,1)\), and
\(\widetilde H_{\mathrm{alg}}(\Ku(X),\Z)\subset L\cap L^{1,1}\)
~\cite{Perry2022IntegralHodge}*{Lemma~5.10}.

\subsection{Freeness and unimodularity}

\begin{prop}\label{prop:L-unimodular}
We have \(K^1_{\mathrm{top}}(X)=0\), and \(L\) is a torsion-free unimodular lattice
of rank \(24\).
\end{prop}

\begin{proof}
The Schubert decomposition of \(G=\Gr(3,10)\), the integral weak Lefschetz theorem,
and Poincar\'e duality imply that \(H^*(X,\Z)\) is free outside the middle degree and vanishes in odd degrees.
In the middle degree, \(H_{19}(X,\Z)=0\) and the universal coefficient theorem give
\(H^{20}(X,\Z)\simeq\Hom(H_{20}(X,\Z),\Z)\), so this group is free as well.
It follows from the Atiyah--Hirzebruch theorem
~\cite{AddingtonThomas2014Hodge}*{Theorem~2.1}
and the vanishing of odd cohomology that \(K^0_{\mathrm{top}}(X)\) is free and
\(K^1_{\mathrm{top}}(X)=0\).

By parts (2) and (3) of the same theorem, for each basis element
\(\alpha\) of \(H^{2k}(X,\Z)\), we may choose a class
\(x_\alpha\in K^0_{\mathrm{top}}(X)\) whose Mukai vector has leading term \(\alpha\).
These classes form a \(\Z\)-basis of \(K^0_{\mathrm{top}}(X)\).
Indeed, the leading term of the Mukai vector of any class is integral,
so successive subtraction of integral linear combinations of the \(x_\alpha\),
in order of increasing degree, eliminates every component.
Linear independence also follows by comparing leading terms.

For \(\alpha\in H^{2k}(X,\Z)\) and \(\beta\in H^{2k'}(X,\Z)\), HRR gives
\[
 \begin{aligned}
 \chi_X(x_\alpha,x_\beta)&=0 &&(k+k'>20),\\
 \chi_X(x_\alpha,x_\beta)&=(-1)^k\int_X\alpha\cup\beta &&(k+k'=20).
 \end{aligned}
\]
Ordering the row basis by increasing degree and the column basis by decreasing degree
gives a block upper triangular Euler matrix whose diagonal blocks are,
up to sign, the matrices of the Poincar\'e pairings.
Ordinary Poincar\'e duality makes each diagonal block unimodular.
Thus the Euler matrix has determinant \(\pm1\), and \(\chi_X\) is perfect.
The Euler matrix for the semiorthogonal decomposition is block upper triangular,
and the exceptional block has determinant one, so the pairing on \(\Ku(X)\) is also unimodular.

On the other hand, the sum of the Betti numbers of \(G\) is 120, and \(b_{20}(G)=b_{22}(G)=10\).
The middle cohomology of \(X\) has rank \(10+22\)
~\cites{DebarreVoisin2010Grassmann,BenedettiSong2021Divisors}.
By weak Lefschetz and duality, the sum of the Betti numbers outside the middle degree is \(120-10-10\), so
\[
 \operatorname{rank}K^0_{\mathrm{top}}(X)=120-10-10+(10+22)=132,\qquad
 \operatorname{rank}L=132-108=24.
\]
\end{proof}

\begin{lem}\label{lem:square-two-index}
Let \(J\) be a unimodular integral lattice, and let \(a\in J\) satisfy \(a^2=2\).
Then \(a\) is primitive, and
\[
 \operatorname{div}_J(a)=1,\qquad
 [J:\Z a\oplus a^\perp]=2,\qquad
 |\operatorname{disc}(a^\perp)|=2.
\]
\end{lem}

\begin{proof}
Since \(a^2=2\), the vector \(a\) is primitive,
and unimodularity implies \((a,J)=\Z\).
The kernel of the surjection \(x\mapsto(a,x)\bmod2\) is \(\Z a\oplus a^\perp\),
so this sublattice has index two.
The discriminant formula for sublattices gives
\[
 2|\operatorname{disc}(a^\perp)|=2^2|\operatorname{disc}J|=4.
\]
\end{proof}

\section{Geometric objects and their Mukai classes}\label{sec:euler-square}

\subsection{The normal bundle and the projected object}

For \(U\in M\), set \(Z_U\colonequals \Gr(3,U)\subset X\),
and let \(i_U\colon Z_U\hookrightarrow X\) denote the inclusion.
Both \(Z_U\) and \(X\) are smooth, of dimensions 9 and 20, respectively.
Thus \(i_U\) is a regular embedding of codimension 11.
Applying adjunction to \(\omega_{Z_U}\simeq\mathcal O_{Z_U}(-6)\) and
\(\omega_X\simeq\mathcal O_X(-9)\) gives
\begin{equation}\label{eq:fiber-normal-determinant}
 \det N_{Z_U/X}
 \simeq\omega_{Z_U}\otimes i_U^*\omega_X^{-1}
 \simeq\mathcal O_{Z_U}(3).
\end{equation}
Using this line bundle, we define
\[
 F_U\colonequals i_{U*}\mathcal O_{Z_U}(3),\qquad E_U\colonequals \delta(F_U)\in\Ku(X),
\]
and write \(v_U\colonequals [E_U]\in L\) for its Mukai class.
Grothendieck duality gives
\[
 F_U^\vee\simeq i_{U*}\mathcal O_{Z_U}[-11],\qquad
 \RHom_X(F_U,A)\simeq R\Gamma(Z_U,i_U^*A)[-11],
\]
where \(A\) is a vector bundle on \(X\) in the second formula.
The twist by \(3\) cancels the factor \(\det N_{Z_U/X}\) arising in duality.

\begin{rem}
This construction has an analogue for lines on a cubic fourfold.
For a smooth cubic fourfold \(Y\subset\mathbf P^5\) and a line \(\ell\subset Y\),
we have \(\omega_Y=\mathcal O_Y(-3)\) and \(\omega_\ell=\mathcal O_\ell(-2)\), hence
\(\det N_{\ell/Y}\simeq\mathcal O_\ell(1)\).
Let \(\delta_Y\colon\Db(Y)\to\Ku(Y)\) be the projection associated with
\(\Db(Y)=\langle\Ku(Y),\mathcal O_Y,\mathcal O_Y(1),\mathcal O_Y(2)\rangle\).
Then the classes of Addington--Thomas can be written as
\[
 \lambda_1=[\delta_Y(\mathcal O_\ell(1))],\qquad
 \lambda_2=[\delta_Y(\mathcal O_\ell(2))]
\]
~\cite{AddingtonThomas2014Hodge}*{Section~2.4}.
Thus both \(\delta_Y(\mathcal O_\ell(1))\) and \(E_U\) are obtained by pushing forward
the determinant of the normal bundle of a subvariety and then projecting to the Kuznetsov component.
The corresponding \(K\)-classes are \(\lambda_1\) and \([E_U]\), respectively.
\end{rem}

\subsection{Computing the Mukai square}

The projected \(K\)-class can be determined from the Euler matrix of the exceptional collection.

\begin{lem}\label{lem:projection-euler}
For \(T,H\in\Db(X)\), set
\[
 G_{ij}\colonequals \chi_X(A_i,A_j),\qquad
 b(T)_i\colonequals \chi_X(A_i,T),\qquad d(T)_i\colonequals \chi_X(T,A_i),
\]
and let \(c(T)\colonequals G^{-1}b(T)\). Then \(c(T)\) is an integral vector, and
\[
 [\delta(T)]=[T]-\sum_i c(T)_i[A_i],
 \qquad
 \chi_X(\delta(T),\delta(H))
 =\chi_X(T,H)-d(T)^{\mathsf t}c(H).
\]
\end{lem}

\begin{proof}
For the projection triangle \(P_T\to T\to\delta(T)\), write
\([P_T]=\sum_i c_i[A_i]\). Semiorthogonality gives
\(\chi_X(A_j,\delta(T))=0\), so
\[
 b(T)_j=\chi_X(A_j,P_T)=\sum_i G_{ji}c_i.
\]
Thus \(Gc=b(T)\). Since \(G\) is an integral upper triangular matrix with all diagonal entries equal to one,
its inverse is also integral, and \(c=c(T)\) is uniquely determined.
The first equality follows from additivity in triangles.

Moreover, since \(P_T\in\mathcal A_X\) and \(\delta(H)\in\mathcal A_X^\perp\),
we have \(\chi_X(P_T,\delta(H))=0\). Hence
\[
 \begin{aligned}
 \chi_X(\delta(T),\delta(H))
 &=\chi_X(T,\delta(H))\\
 &=\chi_X(T,H)-\sum_i c(H)_i\chi_X(T,A_i).
 \end{aligned}
\]
This proves the second equality.
\end{proof}

\begin{lem}\label{lem:F-self-Euler-zero}
For every \(U\in M\), we have
\[
 \chi_X(F_U,F_U)=0.
\]
\end{lem}

\begin{proof}
Since \(Z_U=\Gr(3,U)\subset X\) has codimension 11, GRR gives
\(\ch_k(F_U)=0\) for \(k<11\).
Since \(\ch_k(F_U^\vee)=(-1)^k\ch_k(F_U)\),
every component of the product \(\ch(F_U^\vee)\ch(F_U)\) has codimension at least 22.
As \(22>\dim X=20\), HRR gives
\[
 \chi_X(F_U,F_U)
 =\int_X\ch(F_U^\vee)\ch(F_U)\td(X)=0.\qedhere
\]
\end{proof}

\begin{prop}\label{prop:v-square}
For every \(U\in M\), we have
\[
 (v_U,v_U)=2.
\]
In particular, \(v_U\) is a primitive vector in \(L\), and
\[
 \left|\operatorname{disc}(v_U^\perp)\right|=2.
\]
\end{prop}

\begin{proof}
Applying Lemma~\ref{lem:projection-euler} to \(T=H=F_U\) gives
\[
 \chi_X(E_U,E_U)
 =\chi_X(F_U,F_U)-d(F_U)^{\mathsf t}c(F_U).
\]
The first term vanishes by Lemma~\ref{lem:F-self-Euler-zero}.
With \(Z=\Gr(3,6)\) and \(H_Z=c_1(\mathcal O_Z(1))\), the Euler vectors in the second term are
\[
 \begin{aligned}
 b(F_U)_i&=\int_Z\ch(A_i^\vee|_Z)e^{3H_Z}\td(Z),\\
 d(F_U)_i&=-\int_Z\ch(A_i|_Z)\td(Z).
 \end{aligned}
\]
The first formula follows from pushforward along \(i_U\) and HRR.
The second follows from Grothendieck duality, since the codimension is 11 and
\(\det N_{Z/X}\simeq\mathcal O_Z(3)\).
The formulas, code, and output in Appendix~\ref{sec:computation} give
\[
 d(F_U)^{\mathsf t}G^{-1}b(F_U)=2.
\]
Thus \(\chi_X(E_U,E_U)=-2\), or equivalently, \(v_U^2=2\).

Proposition~\ref{prop:L-unimodular} and Lemma~\ref{lem:square-two-index} imply
that \(v_U\) is primitive and \(|\operatorname{disc}(v_U^\perp)|=2\).
\end{proof}

\section{A relative family and the Mukai homomorphism}\label{sec:incidence-kernel}

We assemble the objects \(E_U\) constructed in the preceding section into a perfect family on \(X\times M\),
and define their common Mukai class \(v\).
Using the first Chern class of the determinant associated with the \(K\)-class of this family,
we define the Mukai homomorphism \(\Theta_{\cE}\colon L\to H^2(M,\Z)\).
The main theorem concerns its restriction to \(v^\perp\).

\subsection{The incidence family}\label{subsec:incidence-family}

Let \(p\colon X\times M\to X\) and \(q\colon X\times M\to M\) be the projections.
The incidence variety is
\[
 \cZ=\Gr_M(3,\mathcal U_M)
 =\{(W,U)\in X\times M\mid W\subset U\}.
\]
For each morphism \(g\colon T\to M\), this relative Grassmannian
parametrizes rank-\(3\) subbundles \(\mathcal W\subset g^*\mathcal U_M\)
with locally free quotient. If \(\sigma|_U=0\), then \(\sigma|_W=0\)
for every \(W\subset U\), which gives the incidence description above.
Write \(i\colon\cZ\hookrightarrow X\times M\) for the closed embedding,
\(p_{\cZ}\colonequals p\circ i\), and \(q_{\cZ}\colonequals q\circ i\).
Then \(q_{\cZ}\) is a smooth morphism of relative dimension nine.
Define
\[
 \mathcal O_{\cZ}(1)\colonequals i^*p^*\mathcal O_X(1),\qquad
 \mathcal O_{\cZ}(3)\colonequals\mathcal O_{\cZ}(1)^{\otimes3},
\]
and set
\[
 \cF\colonequals i_*\mathcal O_{\cZ}(3)
 \simeq i_*\mathcal O_{\cZ}\otimes p^*\mathcal O_X(3),
\]
where the isomorphism follows from the projection formula.
Since \(\cZ\to M\) is flat and \(\mathcal O_{\cZ}(3)\) is a line bundle,
\(\cF\) is also flat over \(M\).
Define \(j_U\colon X\hookrightarrow X\times M\) by \(j_U(x)=(x,U)\),
and let \(k_U\colon Z_U\hookrightarrow\cZ\) be the inclusion of the fiber.

The following diagram is Cartesian:
\[
\begin{tikzcd}[column sep=large, row sep=large]
 Z_U \arrow[r, "k_U"] \arrow[d, hook, "i_U"']
   & \cZ \arrow[d, hook, "i"] \\
 X \arrow[r, hook, "j_U"'] & X\times M
\end{tikzcd}
\]
Since \(q_{\cZ},q\) are flat, this square is Tor-independent.
Base change
~\cite{StacksProject}*{\href{https://stacks.math.columbia.edu/tag/08IB}{Tag~08IB}}
therefore gives
\[
 \begin{aligned}
 Lj_U^*\cF
 &\simeq Lj_U^*Ri_*\mathcal O_{\cZ}(3)\\
 &\simeq Ri_{U*}Lk_U^*\mathcal O_{\cZ}(3)\\
 &\simeq i_{U*}\mathcal O_{Z_U}(3)=F_U.
 \end{aligned}
\]
Here \(Ri_*=i_*\) and \(Ri_{U*}=i_{U*}\), since \(i,i_U\) are closed embeddings,
and \(Lk_U^*\mathcal O_{\cZ}(3)\simeq\mathcal O_{Z_U}(3)\),
since \(\mathcal O_{\cZ}(3)\) is a line bundle.

Moreover, \(\cF\) is perfect because \(X\times M\) is smooth.

\subsection{Relative projection and fibers}

For an admissible subcategory \(\mathcal C\subset\Db(X)\), let \(\mathcal C_M\) be
the idempotent-complete triangulated subcategory generated by
\(p^*C\otimes q^*D\), \(C\in\mathcal C,\ D\in\Db(M)\).
Write \(\mathcal A_{X,M}\colonequals (\mathcal A_X)_M\).
Kuznetsov's faithful base change theorem gives
\begin{equation}\label{eq:relative-sod}
 \Db(X\times M)=
 \langle\Ku(X)_M,\langle A_1\rangle_M,\ldots,\langle A_{108}\rangle_M\rangle
 =\langle\Ku(X)_M,\mathcal A_{X,M}\rangle
\end{equation}
~\cite{Kuznetsov2011BaseChangeSOD}*{Proposition~5.1, Theorem~5.6, Corollary~5.9}.
This base change is faithful because the base is \(\operatorname{Spec}\C\).
Moreover, since \(X\) is smooth and proper, the projections onto the admissible
components have finite cohomological amplitude.

\begin{defn}\label{def:relative-projection}
For \(G\in\operatorname{Perf}(X\times M)\), set
\[
 V_{j,M}(G)\colonequals Rq_*R\mathcal Hom(p^*A_j,G)
\]
and form the evaluation triangle
\begin{equation}\label{eq:relative-mutation}
 q^*V_{j,M}(G)\otimes p^*A_j\longrightarrow G\longrightarrow\mathbb L^M_{A_j}(G).
\end{equation}
Let \(\delta_M(G)\) be the object obtained by iterating these mutations in the order \(j=108,\ldots,1\).
This functor is the \(M\)-linear projection onto the first component of \eqref{eq:relative-sod},
and there is a functorial triangle
\begin{equation}\label{eq:relative-projection-triangle}
 P_G^M\longrightarrow G\longrightarrow\delta_M(G),\qquad
 P_G^M\in\mathcal A_{X,M}.
\end{equation}
\end{defn}

Let \(T\) be a Noetherian scheme of finite type over \(\C\),
let \(g\colon T\to M\) be a morphism of finite Tor-dimension,
and set \(j_g=\operatorname{id}_X\times g\).
Define \(\delta_T\) by the same relative mutations on \(X\times T\).
Since perfect base change is compatible with the evaluation morphism,
pulling back each mutation triangle gives
\[
 Lj_g^*\delta_M(G)\simeq\delta_T(Lj_g^*G).
\]
Since \(q\) and its base change are smooth and proper, pushforward along them
preserves perfectness. Hence each mutation, and thus \(\delta_M,\delta_T\), also preserves perfectness.

\begin{defn}\label{def:projected-incidence-family}
For the incidence family \(\cF\), define
\[
 \cE\colonequals \delta_M(\cF)
\]
and write its projection triangle as
\begin{equation}\label{eq:incidence-projection-triangle}
 \cP\longrightarrow\cF\longrightarrow\cE,\qquad\cP\in\mathcal A_{X,M}.
\end{equation}
\end{defn}

\begin{prop}\label{prop:projected-family-fiber}
The family \(\cE\) is perfect, and for \(U\in M\) and \(j_U\colon X\hookrightarrow X\times M\),
we have
\[
 Lj_U^*\cF\simeq F_U=i_{U*}\mathcal O_{\Gr(3,U)}(3),\qquad
 Lj_U^*\cE\simeq E_U=\delta(F_U).
\]
Moreover, \(v_U=[E_U]\in L\) is independent of \(U\).
\end{prop}

\begin{proof}
The first isomorphism was established in Section~\ref{subsec:incidence-family}.
Base change for the relative projection gives
\(Lj_U^*\cE\simeq\delta(Lj_U^*\cF)\simeq\delta(F_U)=E_U\).
Since \(\cF\) is perfect and \(\delta_M\) preserves perfectness,
\(\cE\) is perfect as well.

Thus \([\cE]\in K^0_{\mathrm{top}}(X\times M)\) is defined,
and \(v_U=j_U^*[\cE]\).
Since \(M\) is a connected complex manifold, any two points \(U_0,U_1\in M\)
can be joined by a path \(\gamma\colon[0,1]\to M\).
The map \((x,t)\mapsto(x,\gamma(t))\) is a homotopy between \(j_{U_0}\) and \(j_{U_1}\).
Homotopy invariance of topological \(K\)-theory yields
\[
 v_{U_0}=j_{U_0}^*[\cE]=j_{U_1}^*[\cE]=v_{U_1}.
\]
\end{proof}

We denote this common class by
\[
 v\colonequals v_U=[E_U]\in L.
\]
By Proposition~\ref{prop:v-square}, we have \(v^2=2\).

\subsection{Basic properties of the Mukai homomorphism}

\begin{defn}\label{def:theta}
For a perfect kernel \(\mathcal K\in\operatorname{Perf}(X\times M)\)
and \(x\in K^0_{\mathrm{top}}(X)\), define
\[
 \begin{aligned}
 \Phi^K_{\mathcal K}(x)&\colonequals q_!([\mathcal K]^\vee\otimes p^*x)
                  \in K^0_{\mathrm{top}}(M),\\
 \Theta_{\mathcal K}(x)&\colonequals c_1(\Phi^K_{\mathcal K}(x))\in H^2(M,\Z).
 \end{aligned}
\]
We apply this construction to \(\mathcal K=\cE,\cF\);
in the main theorem, we restrict \(\Theta_{\cE}\) to \(v^\perp\).
\end{defn}
Defining the determinant of a virtual bundle \(W=[D]-[D']\) by
\(\det W=\det D\otimes(\det D')^\vee\), we have
\(c_1(W)=c_1(\det W)\).
Thus the definition above expresses the Mukai homomorphism, defined by the first
Chern class of a determinant, in terms of topological \(K\)-theory.
If \(x\) is an algebraic class, the corresponding determinant line bundle is algebraic as well.

For example, if \(x=[G]\) with \(G\in\Ku(X)\), then
\(\Phi^K_{\cE}(x)\) is the class of the perfect complex
\(Rq_*R\mathcal Hom(\cE,p^*G)\).
The derived fiber of this perfect complex is \(\RHom_X(E_U,G)\),
and the first Chern class of its determinant is \(\Theta_{\cE}([G])\).
For an arbitrary topological \(K\)-class, the formula above still defines an integral class.

\begin{lem}\label{lem:twist-independent}
For \(x\in v^\perp\), the class \(\Theta_{\cE}(x)\) is unchanged under
\(\cE\mapsto\cE\otimes q^*D\), \(D\in\operatorname{Pic}(M)\).
\end{lem}
\begin{proof}
By the projection formula, the index class after twisting is \(D^\vee\otimes\Phi^K_{\cE}(x)\).
Its rank is \(\chi_X(E_U,x)=-(v,x)=0\), so the assertion follows from
\(c_1(D^\vee\otimes W)=c_1(W)-\operatorname{rank}(W)c_1(D)\).
\end{proof}

\begin{prop}\label{prop:theta-Hodge}
The homomorphism \(\Theta_{\cE}\colon L\to H^2(M,\Z)\) is a morphism of integral
Hodge structures of weight two.
\end{prop}
\begin{proof}
It is defined over \(\Z\), since it is the \(c_1\) of a \(K\)-class.
By GRR, it is given by
\begin{equation}\label{eq:theta-Hodge-GRR}
 \Theta_{\cE}(x)
 =\left[q_*\bigl(\ch(\cE^\vee)p^*\ch(x)p^*\td(X)\bigr)\right]_2.
\end{equation}
For \(x\in L^{a,b}\), we have \(\ch_k(x)\in H^{a+k-1,b+k-1}(X)\).
The classes \(\ch_r(\cE^\vee)\) and \(\td_j(X)\) have types \((r,r),(j,j)\), respectively,
and a contribution to degree two requires \(r+k+j=21\).
After pushforward along a morphism of relative dimension \(20\), the type is therefore
\((a+k-1+r+j-20,b+k-1+r+j-20)=(a,b)\).
\end{proof}

\begin{lem}\label{lem:projection-disappears}
For every \(x\in L\), we have \(\Phi^K_{\cE}(x)=\Phi^K_{\cF}(x)\).
\end{lem}
\begin{proof}
Additivity of the relative mutation triangles gives
\[
 [\cP]=\sum_i[p^*A_i\otimes q^*C_i],\qquad C_i\in K^0_{\mathrm{top}}(M).
\]
Here the virtual classes \(C_i\) incorporate the signs arising from shifts.
We take duals, multiply by \(p^*x\), and apply \(q_!\).
By the projection formula, the image of each summand is
\[
 q_!(p^*(A_i^\vee\otimes x)\otimes q^*C_i^\vee)
 =\chi_X(A_i,x)C_i^\vee=0.
\]
Thus the contribution of \(\cP\) vanishes.
\end{proof}

\subsection{The derived dual of the incidence family}

Set
\[
 \mathcal S_3=p_{\cZ}^*\mathcal U_X,\quad
 \mathcal S_6=q_{\cZ}^*\mathcal U_M,\quad
 \mathcal Q_3=\mathcal S_6/\mathcal S_3.
\]
Then \(\mathcal O_{\cZ}(1)=\det\mathcal S_3^\vee\).
The relative tangent bundle of \(\cZ\to M\) is
\(\mathcal S_3^\vee\otimes\mathcal Q_3\).
Since both factors have rank \(3\),
\[
 \det T_{\cZ/M}
 =(\det\mathcal S_3^\vee)^3(\det\mathcal Q_3)^3,\qquad
 \omega_{\cZ/M}
 =\mathcal O_{\cZ}(-6)\otimes q_{\cZ}^*\mathcal O_M(3).
\]
Here we used \(\det\mathcal Q_3=\det\mathcal S_6\otimes(\det\mathcal S_3)^{-1}\).
Since the smooth varieties \(\cZ,X\times M\) have dimensions \(13,24\), respectively,
\(i\) is a regular embedding of codimension \(11\).
Using \(\omega_X=\mathcal O_X(-9)\), relative adjunction gives
\[
 \det N_{\cZ/(X\times M)}
 =\mathcal O_{\cZ}(3)\otimes q_{\cZ}^*\mathcal O_M(3).
\]
Grothendieck duality yields
\[
 \cF^\vee=i_*\bigl(\mathcal O_{\cZ}(-3)\otimes\det N_{\cZ/(X\times M)}\bigr)[-11]
 =i_*q_{\cZ}^*\mathcal O_M(3)[-11].
\]
The sign from the odd shift and the projection formula give
\begin{equation}\label{eq:theta-GRR}
 \Phi^K_{\cF}(x)=-\mathcal O_M(3)\otimes q_{\cZ!}(p_{\cZ}^*x).
\end{equation}
Restricting the determinant of the relative normal bundle to a fiber recovers
the computation in \eqref{eq:fiber-normal-determinant}.

\section{Vanishing cohomology and the incidence correspondence}\label{sec:vanishing}

\subsection{A rational lift of the vanishing cohomology}

Let \(j\colon X\hookrightarrow\Gr(3,V)\) be the inclusion, and define
\[
 H^{20}(X,\Z)_{\mathrm{van}}
 \colonequals \ker\left(j_*\colon H^{20}(X,\Z)\to H^{22}(\Gr(3,V),\Z)\right).
\]
This lattice has rank \(22\) and discriminant \(11\).
The incidence isometry of Benedetti--Song is
\begin{equation}\label{eq:BS-incidence}
 q_{\cZ*}p_{\cZ}^*\colon
 (H^{20}(X,\Z)_{\mathrm{van}}(9),\,\cdot)
 \xrightarrow{\sim}(h^\perp\subset H^2(M,\Z),\,-q_M)
\end{equation}
~\cite{BenedettiSong2021Divisors}*{Theorem~4.6}.
The Tate twist \((9)\) changes the Hodge weight;
it is distinct from the minus sign in the pairing on the target.

Using the isomorphism
\(\mu_X=\ch(-)\td(X)\colon K^0_{\mathrm{top}}(X)_\Q
\to H^{2\ast}(X,\Q)\), set
\begin{equation}\label{eq:mu-pure}
 x_\gamma\colonequals \mu_X^{-1}(\gamma)
\end{equation}
for \(\gamma\in H^{20}(X,\Q)_{\mathrm{van}}\).
Here \(\ch(x_\gamma)\td(X)=\gamma\) is an equality in the full even cohomology,
not merely in its middle-degree component.

For every ambient class \(j^*\alpha\) of positive degree, we have \(\gamma j^*\alpha=0\).
Indeed, let \(r>0\) be the codimension of \(\alpha\). If \(r>10\), the product vanishes for degree reasons.
If \(1\le r\le10\), a class \(\zeta\in H^{20-2r}(X,\Q)\) of complementary degree
can be written as \(j^*\xi\) by weak Lefschetz, so
\[
 \int_X\gamma j^*\alpha\,\zeta
 =\int_{\Gr(3,V)}j_*\gamma\,\alpha\xi=0.
\]
Poincar\'e duality gives \(\gamma j^*\alpha=0\).
Since \(\td(X)\) is an ambient class with constant term \(1\), we have
\[
 \ch(x_\gamma)=v_X(x_\gamma)=\mu_X(x_\gamma)=\gamma.
\]
Moreover, \(\ch(A_i^\vee)=j^*\alpha_i\), so
\[
 \chi_X(A_i,x_\gamma)=\int_Xj^*\alpha_i\gamma
 =\int_{\Gr(3,V)}\alpha_i j_*\gamma=0.
\]
Thus \(x_\gamma\in L_\Q\).
Denote the image of this lift by \(L_{\mathrm{van},\Q}\subset L_\Q\).

\begin{lem}\label{lem:kappa-pairing}
For \(\gamma,\eta\in H^{20}(X,\Q)_{\mathrm{van}}\), we have
\[
 (x_\gamma,x_\eta)=-\int_X\gamma\eta.
\]
Moreover, \(L_{\mathrm{van},\Q}\) is a nondegenerate sub-Hodge structure of rank \(22\),
and the second summand in
\begin{equation}\label{eq:Hodge-splitting}
 L_\Q=L_{\mathrm{van},\Q}\perp L_{\mathrm{amb},\Q},\qquad
 L_{\mathrm{amb},\Q}\colonequals L_{\mathrm{van},\Q}^{\perp}
\end{equation}
is a Tate Hodge structure of rank \(2\).
\end{lem}
\begin{proof}
Since \(\gamma\) has codimension \(10\),
we have \(\ch(x_\gamma^\vee)=(-1)^{10}\gamma=\gamma\).
Since \(\gamma\eta\) has top degree, only the constant term of the Todd class
contributes in HRR, giving the asserted formula for the pairing.
The orthogonal decomposition follows from the nondegeneracy of the intersection form.
By \(v_X(x_\gamma)=\gamma\) and Definition~\ref{def:mukai-lattice},
this lift is a Hodge embedding of \(H^{20}(X,\Q)_{\mathrm{van}}(9)\).
The off-diagonal Hodge components \(H^{11,9},H^{9,11}\) of \(X\) are entirely
contained in the vanishing cohomology, so the remaining \(24-22=2\)-dimensional
orthogonal complement has only type \((1,1)\).
\end{proof}

\subsection{Comparison of the Mukai homomorphism with the incidence map}

\begin{prop}\label{prop:theta-trans-rational}
\[
 \Theta_{\cE}(x_\gamma)=-q_{\cZ*}p_{\cZ}^*\gamma.
\]
Consequently, the restriction
\(\Theta_{\cE,\Q}|_{L_{\mathrm{van},\Q}}\colon L_{\mathrm{van},\Q}\to h^\perp_\Q\)
is a rational Hodge isometry.
\end{prop}
\begin{proof}
Lemma~\ref{lem:projection-disappears} and GRR give
\[
 \Theta_{\cE}(x_\gamma)
 =\left[q_*\bigl(\ch(\cF^\vee)p^*\ch(x_\gamma)p^*\td(X)\bigr)\right]_2
 =\left[q_*(\ch(\cF^\vee)p^*\gamma)\right]_2.
\]
Since \(\cF^\vee=i_*q_{\cZ}^*\mathcal O_M(3)[-11]\),
the first nonzero component of its Chern character is \(-[\cZ]\), in codimension \(11\).
The minus sign comes from the shift \([-11]\).
Since \(q_*\) lowers codimension by \(20\) and \(\gamma\) has codimension \(10\),
only this lowest-degree component contributes to degree two.
Thus
\[
 \Theta_{\cE}(x_\gamma)
 =-q_*([\cZ]p^*\gamma)=-q_{\cZ*}p_{\cZ}^*\gamma.
\]
By \eqref{eq:BS-incidence} and Lemma~\ref{lem:kappa-pairing},
\[
 q_M(\Theta_{\cE}(x_\gamma),\Theta_{\cE}(x_\eta))
 =-\int_X\gamma\eta=(x_\gamma,x_\eta).
\]
Hence Hodge types and the pairing are preserved.
The lift and the isomorphism here are defined over the rationals;
the comparison of the full integral lattices is carried out in Section~\ref{sec:main-proof}.
\end{proof}

\section{The preimage of the Pl\"ucker polarization}\label{sec:algebraic-lattice}

\subsection{Pairings of two auxiliary classes}

To recover the polarization direction, we use two auxiliary bundles with twist 9,
the first twist outside the exceptional range \(0,\ldots,8\). Set
\[
 B_0=\mathcal O_X(9),\quad B_4=\operatorname{Sym}^4\mathcal U_X^\vee(9),
 \qquad s_i=[\delta(B_i)],\quad m_i=(s_i,v),\quad u_i=2s_i-m_iv.
\]
Since \(v^2=2\), we have \((u_i,v)=0\).

\begin{prop}\label{prop:boundary-computation}
\[
 m_0=-510,\qquad m_4=-97695,\qquad
 ((s_i,s_j))_{i,j=0,4}
 =\begin{pmatrix}146088&43960320\\43960320&27395296750\end{pmatrix}.
\]
Consequently,
\[
 ((u_i,u_j))_{i,j=0,4}
 =22\begin{pmatrix}54^2&54\cdot64135\\54\cdot64135&64135^2\end{pmatrix}.
\]
\end{prop}
\begin{proof}
By Lemma~\ref{lem:projection-euler}, computing \(c_i=G^{-1}b(B_i)\) gives
\[
 \begin{aligned}
 (s_i,s_j)&=-\chi_X(B_i,B_j)+d(B_i)^{\mathsf t}c_j,\\
 (s_i,v)&=-\chi_X(B_i,F_U)+d(B_i)^{\mathsf t}c(F_U).
 \end{aligned}
\]
The entries of \(G,b,d\) are computed by integration on the Grassmannian,
using the hyperplane exact sequence and HRR.
In particular, Serre duality and \(\omega_X=\mathcal O_X(-9)\) give
\[
 \chi_X(B_i,A_{\lambda,t})
 =\chi_X(A_{\lambda,t},\operatorname{Sym}^i\mathcal U_X^\vee),
\]
which vanishes for \(t>0\) by semiorthogonality.
The integration formulas and code in Appendix~\ref{sec:computation} yield the five integers in the proposition.
Finally, the identity
\[
 (u_i,u_j)=4(s_i,s_j)-2m_im_j
\]
gives the second matrix.
\end{proof}

\subsection{Coefficients of the Mukai homomorphism}\label{sec:plucker-preimage}

For \(\lambda\in\Lambda\), append trailing zeros so that it has three entries,
and write \(|\lambda|=\lambda_1+\lambda_2+\lambda_3\).
For \(t\ge0\), set
\(A_{\lambda,t}\colonequals \Sigma^\lambda\mathcal U_X^\vee(t)\) and
\(\mu=\lambda+(t,t,t)\).
Since the weight \((\mu_1,\mu_2,\mu_3,0,0,0)\) is dominant,
the relative Borel--Weil theorem on the Grassmannian
~\cite{Weyman2003Cohomology}*{Corollary~4.1.9} gives
\[
 Rq_{\cZ*}p_{\cZ}^*A_{\lambda,t}
 =Rq_{\cZ*}\Sigma^\mu\mathcal S_3^\vee
 =\Sigma^\mu\mathcal U_M^\vee
\]
with vanishing higher direct images.
In characteristic zero, Weyman's Weyl functor \(K_\mu\) agrees with the Schur functor
\(\Sigma^\mu\) used here.
Let \(r_{\lambda,t}=\dim\Sigma^\mu\C^6\).
The weights of the Schur representation have total degree \(|\mu|=|\lambda|+3t\),
and its character is symmetric in the six Chern roots.
Thus the coefficient of each root in \(c_1\) is the same.
Since the sum of these coefficients is \(|\mu|r_{\lambda,t}\), we obtain
\[
 c_1(\Sigma^\mu\mathcal U_M^\vee)
 =\frac{|\lambda|+3t}{6}\,r_{\lambda,t}h.
\]
Taking into account both \(\mathcal O_M(3)\) and the minus sign in~\eqref{eq:theta-GRR}, we obtain
\begin{equation}\label{eq:relative-BW-theta}
 \Theta_{\cF}([A_{\lambda,t}])=-q_{\lambda,t}h,\qquad
 q_{\lambda,t}\colonequals r_{\lambda,t}
 \left(3+\frac{|\lambda|+3t}{6}\right).
\end{equation}

\begin{lem}\label{lem:direct-incidence-vanishes}
We have \(\Theta_{\cF}([F_{U_0}])=0\) for all \(U_0\in M\).
\end{lem}
\begin{proof}
Both \(\ch(F_{U_0})\) and \(\ch(\cF^\vee)\) vanish in codimension less than \(11\).
By GRR for the family index, every term in the integrand has codimension at least \(22\).
Since \(q\) has relative dimension \(20\), applying \(q_*\) lowers codimension by \(20\),
so the resulting class has codimension at least \(2\).
Thus \(c_1\) vanishes rationally.
Since \(H^2(M,\Z)\) is torsion-free, it also vanishes integrally.
\end{proof}

\begin{prop}\label{prop:theta-auxiliary}
\[
 \Theta_{\cE}(v)=-3h,\quad
 \Theta_{\cE}(s_0)=792h,\quad\Theta_{\cE}(s_4)=178610h.
\]
Consequently, \(\Theta_{\cE}(u_0)=54h\) and \(\Theta_{\cE}(u_4)=64135h\).
\end{prop}
\begin{proof}
Set \(\boldsymbol q=(q_a)_{a=1}^{108}\), where \(q_a=q_{\lambda,t}\)
if \(A_a=A_{\lambda,t}\). Fix \(U\in M\) and write \(c_U=c(F_U)\) and \(c_i=c(B_i)\).
The identity \(v=[F_U]-\sum_a c_{U,a}[A_a]\) and the preceding lemma give
\[
 \Theta_{\cF}(v)=(c_U^{\mathsf t}\boldsymbol q)h,\qquad
 \Theta_{\cF}(s_i)=(-q(B_i)+c_i^{\mathsf t}\boldsymbol q)h.
\]
Here \(q(B_i)\) is the coefficient in~\eqref{eq:relative-BW-theta} corresponding to
\(\lambda=(i)\) and \(t=9\).
Using Weyl's dimension formula and the projection coefficients obtained from the Euler matrix, we compute
\[
 \begin{array}{c|r}
 q(B_0)&4294290\\
 q(B_4)&157237080\\
 c_U^{\mathsf t}\boldsymbol q&-3\\
 c_0^{\mathsf t}\boldsymbol q&4295082\\
 c_4^{\mathsf t}\boldsymbol q&157415690
 \end{array}.
\]
The formulas and code used in this computation are given in Appendix~\ref{sec:computation}.
By Lemma~\ref{lem:projection-disappears},
\(\Theta_{\cF}\) and \(\Theta_{\cE}\) agree on the classes \(v,s_i\in L\),
which proves the first three identities.
Moreover,
\[
 \Theta_{\cE}(u_0)/h=2\cdot792-510\cdot3=54,\qquad
 \Theta_{\cE}(u_4)/h=2\cdot178610-97695\cdot3=64135.
\]
\end{proof}

\subsection{The polarization class and the orthogonal decomposition}

\begin{defn}\label{def:U}
Using the B\'ezout identity \(-22566\cdot54+19\cdot64135=1\), define
\[
 \uDV\colonequals -22566u_0+19u_4.
\]
\end{defn}

\begin{prop}\label{prop:U-lattice}
\[
 (\uDV,v)=0,\quad \uDV^2=22,\quad
 \Theta_{\cE}(\uDV)=h,\quad \uDV\equiv v\pmod{2L}.
\]
Moreover, \(s=(v+\uDV)/2\) is an algebraic class, and
\(\Z v+\Z s\) has Gram matrix
\(\left(\begin{smallmatrix}2&1\\1&6\end{smallmatrix}\right)\)
and discriminant \(11\).
\end{prop}
\begin{proof}
Orthogonality follows from \(u_i\in v^\perp\).
The matrix in Proposition~\ref{prop:boundary-computation} and the B\'ezout identity give
\(\uDV^2=22(-22566\cdot54+19\cdot64135)^2=22\).
Applying the same integral linear combination to Proposition~\ref{prop:theta-auxiliary} gives the image \(h\).
On the other hand,
\[
 \uDV=2(-22566s_0+19s_4)-9652455v,
\]
so \(\uDV\equiv v\pmod{2L}\), and
\[
 s=-22566s_0+19s_4-4826227v
\]
is an integral linear combination of classes of objects.
The identities \((v,s)=(v^2+(v,\uDV))/2=1\) and
\(s^2=(v^2+\uDV^2)/4=6\) give the Gram matrix and its discriminant.
\end{proof}

\begin{lem}\label{lem:ambient-generators}
\[
 L_{\mathrm{amb},\Q}=\Q v\perp\Q\uDV,\qquad
 (v^\perp)_\Q=L_{\mathrm{van},\Q}\perp\Q\uDV.
\]
\end{lem}
\begin{proof}
For the lift \(x_\gamma\) of a vanishing class, we have \(\ch(x_\gamma)\td(X)=\gamma\).
Since \(\ch(F_U^\vee)\) has codimension at least \(11\) and
\(\gamma\) has codimension \(10\), we obtain \(\chi_X(F_U,x_\gamma)=0\).
The Chern characters of the duals of \(B_0,B_4,A_i\) are ambient classes.
Thus \(j_*\gamma=0\) also implies that their Euler pairings with \(x_\gamma\) vanish.
By additivity of the projected classes, \(v,s_i,\uDV\) are orthogonal to \(L_{\mathrm{van},\Q}\).

The orthogonal complement in~\eqref{eq:Hodge-splitting} has dimension \(2\),
and \(v,\uDV\) are orthogonal with nonzero squares \(2,22\).
Thus \(v,\uDV\) form a basis of this orthogonal complement.
Using this basis in the orthogonal decomposition gives both identities in the lemma.
\end{proof}

\section{Integral Hodge isometry and Hodge classes}\label{sec:main-proof}

\subsection{Proof of the main theorem}

\begin{proof}[Proof of Theorem~\ref{thm:main}]
Fix \(\sigma\) for which \(X,M\) are smooth.
We first prove that the rationalized map is an isometry, and then compute the index of its image in the integral lattice.

\emph{Step 1: Rational orthogonal decompositions.}
Lemma~\ref{lem:ambient-generators} and \(q_M(h)=22\ne0\) give
\[
 (v^\perp)_\Q=L_{\mathrm{van},\Q}\mathbin{\perp}\Q\uDV,
 \qquad
 H^2(M,\Q)=H^2(M,\Q)_{\mathrm{prim}}\mathbin{\perp}\Q h.
\]
On the first summands, Proposition~\ref{prop:theta-trans-rational} gives a rational Hodge isometry.
On the second summands, Proposition~\ref{prop:U-lattice} gives
\[
 \Theta_{\cE}(\uDV)=h,\qquad \uDV^2=22=q_M(h).
\]
Thus the pairing is preserved on the second summands as well, and the cross terms vanish on both sides.
Since each summand maps isomorphically onto the corresponding orthogonal summand,
\[
 \Theta_{\cE,\Q}:(v^\perp)_\Q\xrightarrow{\sim}H^2(M,\Q)
\]
is a rational Hodge isometry.

\emph{Step 2: Comparison of the integral lattices.}
By Definition~\ref{def:theta},
\(\Theta_{\cE}(v^\perp)\subset H^2(M,\Z)\).
Since \(v^\perp\) is torsion-free and the rationalized map is injective,
the integral map is also injective. Since its rationalization is surjective,
its image is a sublattice of finite index. Denote this index by \(D\).
By Proposition~\ref{prop:v-square}, we have
\(|\operatorname{disc}(v^\perp)|=2\), and Section~\ref{subsec:notation-geometry} gives
\(|\operatorname{disc}H^2(M,\Z)|=2\).
The preservation of the pairing in Step~1 and the discriminant formula for a sublattice give
\[
 2=|\operatorname{disc}(v^\perp)|
  =D^2|\operatorname{disc}H^2(M,\Z)|=2D^2.
\]
Hence \(D=1\), and \(\Theta_{\cE}|_{v^\perp}\) is an integral Hodge isometry.
This integral isometry maps \(\uDV\) to \(h\).
Since \(h\) is primitive of divisibility \(2\), the class \(\uDV\) is primitive
in \(v^\perp\) and satisfies \(\operatorname{div}_{v^\perp}(\uDV)=2\).
\end{proof}

\subsection{The Hodge lattice at a very general point}

Define the lattice of integral Hodge classes and the transcendental lattice by
\[
 \NHdg\colonequals L\cap L^{1,1},\qquad T\colonequals \NHdg^\perp.
\]
The sublattice \(\NHdg\) is saturated.
Indeed, if \(rx\in L^{1,1}\) for a positive integer \(r\),
then \(x\in L^{1,1}\) by linearity of the Hodge decomposition.

\begin{prop}\label{prop:N-lattice}
For a very general \(\sigma\) in the parameter space of trivectors
for which both \(X_3^\sigma\) and \(X_6^\sigma\) are smooth, we have
\[
 \NHdg=\widetilde H_{\mathrm{alg}}(\Ku(X),\Z)=\Z v\oplus\Z s,\qquad
 |\operatorname{disc}\NHdg|=|\operatorname{disc}T|=11.
\]
Moreover, \(\Theta_{\cE}|_T\colon T\xrightarrow{\sim}h^\perp\) is
an integral Hodge isometry.
\end{prop}
\begin{proof}
At a very general point, the vanishing part has no rational Hodge classes
~\cite{BenedettiSong2021Divisors}*{Section~3.1}.
Thus~\eqref{eq:Hodge-splitting} gives
\(\NHdg\otimes\Q=L_{\mathrm{amb},\Q}\) and \(\operatorname{rank}\NHdg=2\).
By Proposition~\ref{prop:U-lattice},
\[
 K\colonequals \Z v+\Z s\subset
 \widetilde H_{\mathrm{alg}}(\Ku(X),\Z)\subset\NHdg,\qquad
 |\operatorname{disc}K|=11.
\]
The index \(e=[\NHdg:K]\) satisfies
\(11=e^2|\operatorname{disc}\NHdg|\).
Since \(|\operatorname{disc}\NHdg|\) is an integer, we have \(e^2\mid11\), and hence \(e=1\).
Both inclusions are therefore equalities.

By saturation and unimodularity, the map
\(L\to\NHdg^\vee\) induced by the pairing is surjective with kernel \(T\).
Indeed, every integral linear form on the primitive subgroup \(\NHdg\) extends to
\(L\), and \(L=L^\vee\).
It follows that
\[
 [L:\NHdg\oplus T]=|\NHdg^\vee/\NHdg|=11.
\]
Taking the discriminant of the orthogonal sum yields
\(11|\operatorname{disc}T|=11^2|\operatorname{disc}L|=121\).
Finally, since \(\uDV=2s-v\), we have
\(T=v^\perp\cap\uDV^\perp\), which the isometry of the main theorem maps onto \(h^\perp\).
\end{proof}

For the same very general \(\sigma\), the integral lattice \(T\oplus\Z\uDV\) satisfies
\[
 [v^\perp:T\oplus\Z\uDV]^2
 =\frac{11\cdot22}{2}=121,
\]
so its index is \(11\).

\section{Comparison with a twisted K3 surface at
\texorpdfstring{discriminant \(24\)}{discriminant 24}}\label{sec:D24}

\subsection{The twisted K3 surface and its Mukai lattice}

Let \([\sigma]\) be a very general point of \(\mathcal D^{1,6,10}\) as defined in the Introduction,
and assume that \(X,M\) are smooth.
Benedetti--Song construct a degree-\(6\) K3 surface \((S,\ell)\),
a class \(\beta\in\operatorname{Br}(S)\) of order \(2\), and a lift \(B_S\)
~\cite{BenedettiSong2021Divisors}*{Section~5.3, Lemma~5.15}.
These satisfy
\[
 \operatorname{Pic}(S)=\Z\ell,\quad \ell^2=6,\quad
 2B_S\in H^2(S,\Z),\quad B_S^2=B_S\cdot\ell=\frac12.
\]

Here \(\operatorname{Br}(S)=H^2(S,\mathcal O_S^\times)_{\mathrm{tors}}\) is
the Brauer group defined using analytic cohomology.
A \(\beta\)-twisted sheaf consists of local sheaves with gluing isomorphisms
whose triple products are scalar multiples of the identity given by a 2-cocycle representing \(\beta\).
A rational class \(B_S\) is a \(B\)-field lift if
\(\exp(B_S^{0,2})=\beta\) under the map arising from the exponential sequence
~\cite{HuybrechtsStellari2005TwistedK3}*{Proposition~1.2}.
We fix this lift throughout.

\begin{defn}\label{def:twisted-Mukai}
The underlying group and pairing of the twisted Mukai lattice are
\[
 \widetilde H(S,B_S,\Z)=H^0(S,\Z)\oplus H^2(S,\Z)\oplus H^4(S,\Z),\qquad
 ((r,c,a),(r',c',a'))=c\cdot c'-ra'-r'a.
\]
For \(0\ne\omega_S\in H^{2,0}(S)\), its \((2,0)\)-line is
\(\C(0,\omega_S,\omega_S B_S)\), its \((0,2)\)-line is the complex conjugate,
and its \((1,1)\)-part is their orthogonal complement.
\end{defn}
This is the convention of Huybrechts--Stellari and Benedetti--Song, using
\(e^{B_S}(0,\omega_S,0)=(0,\omega_S,\omega_S B_S)\)
~\cite{HuybrechtsStellari2005TwistedK3}*{Definition~2.3},
~\cite{BenedettiSong2021Divisors}*{Section~5.3}.
The map \(e^{B_S}\) is a rational isometry, but need not be an integral automorphism
of the fixed lattice \(H^*(S,\Z)\).
The underlying integral lattice is the usual Mukai lattice
\(\Hyp^{\oplus4}\oplus E_8(-1)^{\oplus2}\).

Set
\[
 w=(2,2B_S,0).
\]
Since \(2B_S\) is integral, so is \(w\), and
\[
 (w,(0,\omega_S,\omega_S B_S))=0,\qquad
 w^2=4B_S^2=2,\qquad (w,(0,\ell,0))=1.
\]
Thus \(w\) is a primitive class of type \((1,1)\).

\subsection{The moduli description and the two Mukai homomorphisms}

Benedetti--Song give an isomorphism
\[
 \iota\colon
 M\xrightarrow{\sim}M_{S,B_S,\ell}(w)
\]
~\cite{BenedettiSong2021Divisors}*{Theorem~5.19}.
The right-hand side is the moduli space of \(\beta\)-twisted \(\ell\)-stable sheaves of class \(w\).
There is no rank-\(1\) torsion-free twisted sheaf for nonzero \(\beta\).
Indeed, its double dual would be a twisted line bundle,
whose transition functions would express \(\beta\) as a coboundary.
Consequently, a rank-\(2\) torsion-free twisted sheaf has no subsheaf of smaller positive rank.
Such sheaves are slope-stable, hence Gieseker-stable.
It follows that \(\ell\) is \(w\)-generic.

To apply Yoshioka's theorem, choose a Brauer--Severi \(\mathbb P^1\)-bundle
\(p\colon Y\to S\) whose class \(w(Y)\in H^2(S,\mu_2)\) is
\([2B_S\bmod 2]\), and take the auxiliary bundle \(G=\epsilon(Y)\)
~\cite{Yoshioka2006TwistedSheaves}*{Theorem~3.14 and Section~3.3}.
The image of \(w(Y)\) in the Brauer group is \(\beta\).
For the \(Y\)-sheaf \(E_Y\) corresponding to a \(\beta\)-twisted sheaf \(E\),
the Mukai vectors satisfy \(v^{B_S}(E)=e^{B_S}v_G(E_Y)\)
~\cite{Yoshioka2006TwistedSheaves}*{Lemma~3.3 and Remark~3.2}.
Thus the moduli problem of class \(w\) corresponds to the vector
\(w_Y=e^{-B_S}w=(2,0,-1/2)\) in Yoshioka's notation.
Taking \(r=2,\ \xi=2B_S,\ a=0\) gives \(v_Y=e^{\xi/r}w_Y=w\)
in the integral lattice. The rank, primitivity, square, and genericity
conditions checked above allow us to apply
~\cite{Yoshioka2006TwistedSheaves}*{Theorem~3.19(i)} to \(v_Y^\perp=w^\perp\).

Let \(\theta_w\) be the map defined using a twisted quasi-universal family
and normalized by its similitude, and set
\[
 \theta_S\colonequals \iota^*\theta_w\colon
 w^\perp\xrightarrow{\sim}H^2(M,\Z).
\]
This is an integral Hodge isometry.
For \(\theta_w\), we follow Yoshioka's conventions for the dual, \(B\)-field,
Todd class, and similitude~\cite{Yoshioka2006TwistedSheaves}*{Section~3.3, p.~23}.
The map on \(w^\perp\) is independent of the normalization of the quasi-universal family.

Composing the isometry of Theorem~\ref{thm:main} with \(\theta_S^{-1}\), we obtain
\begin{equation}\label{eq:D24-orthogonal-isometry}
 f\colonequals \theta_S^{-1}\circ\Theta_{\cE}|_{v^\perp}\colon
 v^\perp\xrightarrow{\sim}w^\perp.
\end{equation}
This is an integral Hodge isometry through the common lattice \(H^2(M,\Z)\).
Both orthogonal complements have rank \(23\).
We extend \(f\) to the full rank-\(24\) Mukai lattices on both sides by setting \(v\mapsto\pm w\).

\subsection{Extension from the orthogonal complements to the full lattices}

\begin{lem}\label{lem:square-two-gluing}
Let \(L,L'\) be unimodular lattices, and let
\(a\in L\) and \(a'\in L'\) be primitive vectors of square 2.
Every isometry
\[
 g\colon a^\perp\xrightarrow{\sim}(a')^\perp
\]
extends to two isometries
\[
 \widetilde g_\pm\colon L\xrightarrow{\sim}L',
\qquad
 \widetilde g_\pm(a)=\pm a'.
\]

If, moreover, \(L,L'\) carry weight-two Hodge structures compatible with the pairings,
\(a,a'\) are \((1,1)\)-classes, and \(g\) is a Hodge isometry,
then \(\widetilde g_\pm\) are also Hodge isometries.
\end{lem}

\begin{proof}
Set \(K=a^\perp\) and \(K'=(a')^\perp\).
By Lemma~\ref{lem:square-two-index},
\([L:\Z a\oplus K]=2\) and \(|K^\vee/K|=2\); the same equalities hold for \(L'\).
Choose \(z\in L\) with \((z,a)=1\), and write
\[
 z=\frac a2+k^\vee,\qquad k^\vee=z-\frac a2\in K_\Q.
\]
For \(y\in K\), we have \((k^\vee,y)=(z,y)\in\Z\), so
\(k^\vee\in K^\vee\).
If \(k^\vee\in K\), then \(a/2\in L\), contrary to primitivity.
Thus \([k^\vee]\) is the unique nonzero element of the discriminant group of order 2, and
\[
 L=(\Z a\oplus K)+\Z\left(\frac a2+k^\vee\right).
\]
Similarly, we can write
\(L'=(\Z a'\oplus K')+\Z(a'/2+(k')^\vee)\).

Since \(g_\Q\) maps \(K^\vee\) isomorphically onto \((K')^\vee\),
it identifies the unique nonzero elements of the discriminant groups, giving
\(g_\Q(k^\vee)-(k')^\vee\in K'\).
Define maps on the rational orthogonal sums by
\[
 \widetilde g_\pm^\Q(ra+x)=\pm ra'+g_\Q(x).
\]
These preserve the pairings.
The generator adjoined to the sublattice above maps to
\[
 \pm\frac{a'}2+g_\Q(k^\vee)
 \equiv\frac{a'}2+(k')^\vee\pmod{\Z a'\oplus K'}.
\]
This congruence also holds for the minus sign, since \(-a'/2-a'/2=-a'\) belongs to \(\Z a'\).
Thus both maps send \(L\) into \(L'\).
Applying the same argument to their inverses proves surjectivity.
Since both \(a\mapsto\pm a'\) and \(g\) preserve Hodge type,
the rational orthogonal decompositions show that the extensions also preserve the Hodge structures.
\end{proof}

\begin{proof}[Proof of Theorem~\ref{thm:D24-main}]
Set \(L=\widetilde H(\Ku(X),\Z)\) and \(L'=\widetilde H(S,B_S,\Z)\).
By Proposition~\ref{prop:L-unimodular} and the description of the K3 Mukai lattice,
both lattices are unimodular.
The vectors \(v,w\) are primitive of type \((1,1)\), with \(v^2=w^2=2\),
and~\eqref{eq:D24-orthogonal-isometry} gives a Hodge isometry of their orthogonal complements.
Applying Lemma~\ref{lem:square-two-gluing} with \(a=v,\ a'=w,\ g=f\), we obtain
\[
 \varphi_\pm\colon L\xrightarrow{\sim}L',\qquad
 \varphi_\pm(v)=\pm w,\qquad
 \varphi_\pm|_{v^\perp}=f.
\]
\end{proof}

A Fourier--Mukai equivalence \(\Ku(X)\simeq\Db(S,\beta)\) preserves the Euler pairing
and Hodge structure on topological \(K\)-theory, and hence induces an integral Hodge isometry
of Mukai lattices
~\cite{Perry2022IntegralHodge}*{Section~5.1, Definition~6.4}.
Theorem~\ref{thm:D24-main} provides the isomorphism of Hodge lattices required for such an equivalence,
in a form compatible with \(f\) and \(v\mapsto\pm w\).

\appendix
\section{Macaulay2 calculations}\label{sec:computation}

The numerical values used in Propositions~\ref{prop:v-square}, \ref{prop:boundary-computation},
and \ref{prop:theta-auxiliary} can be computed with the following code.
The inputs are the twelve partitions in \(\Lambda\), the twists \(0,\ldots,8\),
and \(F_U,B_0,B_4\).
Integration using Schubert2 and linear algebra over \(\Q\) yield
the pairings and the values of the Mukai homomorphism.

\subsection{The Euler matrix and three projections}

Let \(G_{\mathrm{amb}}\colonequals \Gr(3,10)\) and \(Z\colonequals \Gr(3,6)\),
and denote their respective Pl\"ucker classes by \(H,H_Z\).
Write \(A_a=A_{\lambda,t}\), \(1\le a\le108\),
in the order of the exceptional collection.
The hyperplane exact sequence and HRR give
\[
 G_{ab}\colonequals \chi_X(A_a,A_b)
 =\int_{G_{\mathrm{amb}}}
 \ch(A_a^\vee)\ch(A_b)(1-e^{-H})\td(G_{\mathrm{amb}}),
\]
where the integrand uses the Chern characters of the corresponding Schur bundles on the ambient Grassmannian.
The entries with \(b<a\) vanish by semiorthogonality.
The entries with \(a\le b\), including the diagonal entries, are computed by the formula above.
For entries with the same difference of twists, we use the common values
\[
 G_{\mathrm{diff}}(d)_{\lambda,\mu}
 \colonequals \chi_X(A_{\lambda,0},A_{\mu,0}(d)).
\]
Computing \(G\) requires \(d=0,\ldots,8\).

We use the notation \(b(T),d(T),c(T)\) from Lemma~\ref{lem:projection-euler}.
For \(T=F_U,B_0,B_4\), we compute the Euler characteristics below
and substitute them into the pairing formula of that lemma.

For \(F_U=i_{U*}\mathcal O_Z(3)\), the two required Euler vectors are
\[
 \begin{aligned}
 b(F_U)_{\lambda,t}
 &=\int_Z\ch(A_{\lambda,t}^\vee|_Z)e^{3H_Z}\td(Z),\\
 d(F_U)_{\lambda,t}\colonequals \chi_X(F_U,A_{\lambda,t})
 &=-\int_Z\ch(A_{\lambda,t}|_Z)\td(Z).
 \end{aligned}
\]
By Lemmas~\ref{lem:projection-euler} and \ref{lem:F-self-Euler-zero},
\[
 \chi_X(E_U,E_U)=-c(F_U)^{\mathsf t}d(F_U)=-2.
\]

\subsection{Pairings of the two auxiliary classes}

For \(i\in\{0,4\}\), set
\(B_i=\operatorname{Sym}^i\mathcal U_X^\vee(9)\) and
\(s_i=[\delta(B_i)]\).
Writing \(\mu_i=(i,0,0)\), we have
\[
 b(B_i)_{\lambda,t}=G_{\mathrm{diff}}(9-t)_{\lambda,\mu_i}.
\]
For \(t=0\), we additionally need entries with twist difference nine.
The code computes only the two additional columns corresponding to \(\mu_0,\mu_4\).
By Serre duality and semiorthogonality, the reverse pairing is
\[
 d(B_i)_{\lambda,t}\colonequals \chi_X(B_i,A_{\lambda,t})
 =\begin{cases}
 G_{\mathrm{diff}}(0)_{\lambda,\mu_i}&t=0,\\
 0&t>0
 \end{cases}.
\]
No further integration is required. We also compute the two values
\[
 e_i\colonequals \chi_X(B_i,F_U)
 =\int_Z\ch(B_i^\vee|_Z)e^{3H_Z}\td(Z).
\]
Since twisting both arguments by the same line bundle leaves the Euler pairing unchanged,
\(\chi_X(B_i,B_j)=G_{\mathrm{diff}}(0)_{\mu_i,\mu_j}\).
Thus the pairings after projection are
\[
 \begin{aligned}
 (s_i,v)&=-e_i+d(B_i)^{\mathsf t}c(F_U),\\
 (s_i,s_j)&=-G_{\mathrm{diff}}(0)_{\mu_i,\mu_j}
               +d(B_i)^{\mathsf t}c(B_j).
 \end{aligned}
\]
By symmetry of the Mukai pairing, it suffices to compute the five values
\[
 (s_0,v),(s_4,v),(s_0,s_0),(s_0,s_4),(s_4,s_4).
\]
For \(m_i=(s_i,v)\) and \(u_i=2s_i-m_iv\), we have
\[
 (u_i,u_j)=4(s_i,s_j)-2m_im_j.
\]

\subsection{Coefficients of the Mukai homomorphism}

Set \(\nu=(\lambda_1+t,\lambda_2+t,\lambda_3+t,0,0,0)\).
By Weyl's dimension formula, the coefficients in~\eqref{eq:relative-BW-theta} are
\[
 r_{\lambda,t}=\prod_{1\le a<b\le6}
 \frac{\nu_a-\nu_b+b-a}{b-a},\qquad
 q_{\lambda,t}=r_{\lambda,t}
 \left(3+\frac{|\lambda|+3t}{6}\right).
\]
For the column vector \(\boldsymbol q=(q_{\lambda,t})\),
Lemmas~\ref{lem:direct-incidence-vanishes} and \ref{lem:projection-disappears},
together with the formula for projected classes, give
\[
 \begin{aligned}
 \Theta_{\cE}(v)/h&=c(F_U)^{\mathsf t}\boldsymbol q,\\
 \Theta_{\cE}(s_i)/h&=-q_{\mu_i,9}+c(B_i)^{\mathsf t}\boldsymbol q.
 \end{aligned}
\]
Finally, we substitute \(u_i=2s_i-m_iv\) and
\(\uDV=-22566u_0+19u_4\).

The main variables in the formulas and the code correspond as follows.
Since the code indexes \(\Lambda\) starting at zero,
the partitions \((0)\) and \((4)\) occupy positions 0 and 9, respectively.
\[
\begin{array}{c|c}
\text{Formula}&\text{Code}\\ \hline
G,\ G_{\mathrm{diff}}&\texttt{Gmat},\ \texttt{Gdiff}\\
b(F_U),d(F_U),c(F_U)&\texttt{bF},\ \texttt{dF},\ \texttt{cF}\\
c(B_0),c(B_4)&\texttt{c0},\ \texttt{c4}\\
d(B_0),d(B_4)&\texttt{d0},\ \texttt{d4}\\
\boldsymbol q&\texttt{qExceptional}\\
(s_0,v),(s_4,v)&\texttt{m0},\ \texttt{m4}\\
(s_0,s_0),(s_0,s_4),(s_4,s_4)&\texttt{s00},\ \texttt{s04},\ \texttt{s44}\\
\Theta_{\cE}(s_0)/h,\ \Theta_{\cE}(s_4)/h
 &\texttt{thetaS0},\ \texttt{thetaS4}\\
\Theta_{\cE}(u_0)/h,\ \Theta_{\cE}(u_4)/h
 &\texttt{thetaU0},\ \texttt{thetaU4}
\end{array}
\]

\subsection{Running the code}
\raggedbottom

The calculations were performed with Macaulay2 1.26.06 and Schubert2 0.7.

The vanishing entries dictated by semiorthogonality, \(\chi_X(F_U,F_U)=0\), and
\(\Theta_{\cF}([F_U])=0\) are inputs proved in the main text.

\begin{lstlisting}[
 basicstyle=\ttfamily\footnotesize,
 lineskip=-0.6pt,
 caption={Euler pairing and Mukai homomorphism},
 label={lst:dv-integral-hodge}
]
loadPackage "Schubert2";

(
chLinePower = (L,t) -> (
    if t >= 0 then (ch L)^t else (ch dual L)^(-t)
);

-- G=Gr(3,10), Z=Gr(3,6), with their Plucker line bundles.
G = flagBundle({7},10);
UG = (bundles G)#0;
LG = dual exteriorPower_3 UG;
Z = flagBundle({3},6);
UZ = (bundles Z)#0;
LZ = dual exteriorPower_3 UZ;

-- HRR on the hyperplane X, using [O_X]=[O_G]-[O_G(-1)].
chiX = chB -> integral(chB * (1-chLinePower(LG,-1)) * todd(G));

Lambda = {
    {0}, {1}, {1,1}, {2}, {2,1}, {2,2},
    {3}, {3,1}, {3,2}, {4}, {4,1}, {4,2}
};
nLambda = #Lambda;
nBlocks = 9;
nExceptional = nLambda*nBlocks;

chBaseG = apply(Lambda, lam -> ch(schur_lam(dual UG)));
chDualBaseG = apply(Lambda, lam -> ch(schur_lam UG));
chBaseZ = apply(Lambda, lam -> ch(schur_lam(dual UZ)));
chDualBaseZ = apply(Lambda, lam -> ch(schur_lam UZ));

-- (Gdiff#d)_(a,b)=chi_X(A_(a,0),A_(b,d)).
Gdiff = apply(nBlocks, d -> matrix(QQ, apply(nLambda, a ->
    apply(nLambda, b -> chiX(
        chDualBaseG#a * chBaseG#b * chLinePower(LG,d)
    )))));
G0 = Gdiff#0;

-- Index i=nLambda*t+a; semiorthogonality gives the lower zeros.
Gmat = matrix(QQ, apply(nExceptional, i -> apply(nExceptional, j ->
    if j < i then 0_QQ else
    (Gdiff#(j//nLambda-i//nLambda))_(i%nLambda,j%nLambda))));

-- F=i_*O_Z(3): b_i=chi(A_i,F), d_i=chi(F,A_i).
bF = transpose matrix(QQ, {apply(nExceptional, i -> integral(
    chDualBaseZ#(i%nLambda) * chLinePower(LZ,3-i//nLambda)
    * todd(Z)))});
dF = transpose matrix(QQ, {apply(nExceptional, i -> -integral(
    chBaseZ#(i%nLambda) * chLinePower(LZ,i//nLambda) * todd(Z)))});
cF = solve(Gmat,bF);

-- chi(F,F)=0 by GRR and codim_X Z=11.
chiEE = -(transpose(cF)*dF)_(0,0);
vSquare = -chiEE;

-- Weyl dimension formula for a GL(6)-representation.
gl6Dimension = mu -> (
    nu := take(mu | {0,0,0,0,0,0},6);
    value := 1_QQ;
    for i from 0 to 4 do (
        for j from i+1 to 5 do (
            value = value * (nu#i-nu#j+j-i)/(j-i);
        );
    );
    value
);

-- Coefficient q_(lambda,t) in Theta_F(A_(lambda,t))=-q_(lambda,t)h.
relativeBWCoefficient = (lambda,t) -> (
    mu := apply(take(lambda | {0,0,0},3), x -> x+t);
    gl6Dimension(mu) * (3_QQ + (sum mu)*1_QQ/6)
);
qExceptional = transpose matrix(QQ, {apply(nExceptional, i ->
    relativeBWCoefficient(Lambda#(i%nLambda),i//nLambda))});

-- Theta_F([F_U])=0, so projection gives Theta(v).
thetaV = (transpose(cF)*qExceptional)_(0,0);

-- Auxiliary bundles B_0=O_X(9), B_4=Sym^4(U_X^vee)(9).
auxiliaryBVector = k -> transpose matrix(QQ, {apply(nExceptional, i -> (
    a := i%nLambda;
    t := i//nLambda;
    if t == 0 then
        chiX(chDualBaseG#a * chBaseG#k * chLinePower(LG,9))
    else (Gdiff#(9-t))_(a,k)
))});

-- Serre duality and semiorthogonality give d_i=0 for t>0.
auxiliaryDVector = k -> transpose matrix(QQ, {apply(nExceptional, i ->
    if i < nLambda then G0_(i,k) else 0_QQ)});
chiAuxiliaryF = k -> integral(
    chDualBaseZ#k * chLinePower(LZ,3-9) * todd(Z)
);

index0 = 0;
index4 = 9; -- Positions of {0} and {4} in Lambda.
b0 = auxiliaryBVector(index0);
b4 = auxiliaryBVector(index4);
c0 = solve(Gmat,b0);
c4 = solve(Gmat,b4);
d0 = auxiliaryDVector(index0);
d4 = auxiliaryDVector(index4);

-- Mukai pairings after projection.
m0 = -(chiAuxiliaryF(index0) - (transpose(d0)*cF)_(0,0));
m4 = -(chiAuxiliaryF(index4) - (transpose(d4)*cF)_(0,0));
s00 = -(G0_(index0,index0) - (transpose(d0)*c0)_(0,0));
s04 = -(G0_(index0,index4) - (transpose(d0)*c4)_(0,0));
s44 = -(G0_(index4,index4) - (transpose(d4)*c4)_(0,0));

qB0 = relativeBWCoefficient(Lambda#index0,9);
qB4 = relativeBWCoefficient(Lambda#index4,9);
sumQ0 = (transpose(c0)*qExceptional)_(0,0);
sumQ4 = (transpose(c4)*qExceptional)_(0,0);
thetaS0 = sumQ0-qB0;
thetaS4 = sumQ4-qB4;

-- u_i=2s_i-m_i v and u_DV=-22566u_0+19u_4.
thetaU0 = 2*thetaS0-m0*thetaV;
thetaU4 = 2*thetaS4-m4*thetaV;
u00 = 4*s00-2*m0^2;
u04 = 4*s04-2*m0*m4;
u44 = 4*s44-2*m4^2;
aBezout = -22566;
bBezout = 19;
thetaUDV = aBezout*thetaU0+bBezout*thetaU4;
vCoefficient = aBezout*m0+bBezout*m4;
uDVSquare = (
    aBezout^2*u00
    + 2*aBezout*bBezout*u04
    + bBezout^2*u44
);

scan({
    {"chi(E,E) = ", chiEE},
    {"(v,v) = ", vSquare},
    {"Theta(v)/h = ", thetaV},
    {"q(B_0) = ", qB0},
    {"c_0^t q = ", sumQ0},
    {"(s_0,v) = ", m0},
    {"(s_0,s_0) = ", s00},
    {"Theta(s_0)/h = ", thetaS0},
    {"Theta(u_0)/h = ", thetaU0},
    {"(u_0,u_0) = ", u00},
    {"q(B_4) = ", qB4},
    {"c_4^t q = ", sumQ4},
    {"(s_4,v) = ", m4},
    {"(s_4,s_4) = ", s44},
    {"Theta(s_4)/h = ", thetaS4},
    {"Theta(u_4)/h = ", thetaU4},
    {"(u_4,u_4) = ", u44},
    {"(s_0,s_4) = ", s04},
    {"(u_0,u_4) = ", u04},
    {"(u_DV,u_DV) = ", uDVSquare},
    {"-22566*(s_0,v) + 19*(s_4,v) = ", vCoefficient},
    {"Theta(u_DV)/h = ", thetaUDV}
}, item -> << item#0 << item#1 << endl);
);
\end{lstlisting}

\begingroup
\renewcommand{\lstlistingname}{Output}
\begin{lstlisting}[
 basicstyle=\ttfamily\footnotesize,
 lineskip=-0.6pt,
 caption={[]Results of Code~\ref{lst:dv-integral-hodge}}
]
chi(E,E) = -2
(v,v) = 2
Theta(v)/h = -3
q(B_0) = 4294290
c_0^t q = 4295082
(s_0,v) = -510
(s_0,s_0) = 146088
Theta(s_0)/h = 792
Theta(u_0)/h = 54
(u_0,u_0) = 64152
q(B_4) = 157237080
c_4^t q = 157415690
(s_4,v) = -97695
(s_4,s_4) = 27395296750
Theta(s_4)/h = 178610
Theta(u_4)/h = 64135
(u_4,u_4) = 90492560950
(s_0,s_4) = 43960320
(u_0,u_4) = 76192380
(u_DV,u_DV) = 22
-22566*(s_0,v) + 19*(s_4,v) = 9652455
Theta(u_DV)/h = 1
\end{lstlisting}
\endgroup

\end{document}